\documentclass{amsart}
\usepackage{amsmath,amsthm,amssymb,amscd,enumerate,tikz}
\usepackage{wrapfig,setspace,mathrsfs,comment,color}
\usepackage{stmaryrd}

\usepackage[all]{xy}

\numberwithin{equation}{section}

\SelectTips{eu}{12}

\newtheorem{thm}{Theorem}[section]

\newtheorem{theorem}[thm]{Theorem}
\newtheorem{proposition}[thm]{Proposition}
\newtheorem{lemma}[thm]{Lemma}
\newtheorem{corollary}[thm]{Corollary}

\newtheorem{mthm}{Theorem}

\theoremstyle{definition}
\newtheorem{definition}[thm]{Definition}
\newtheorem{example}[thm]{Example}

\newtheorem{setting}[thm]{Setting}

\theoremstyle{remark}

\newtheorem{remark}[thm]{Remark}

\newcommand{\QQ}{\mathbb{Q}}
\newcommand{\RR}{\mathbb{R}}
\newcommand{\ZZ}{\mathbb{Z}}
\newcommand{\NN}{\mathbb{N}}

\newcommand{\CC}{\mathcal{C}}
\newcommand{\hh}{\mathrm{h}}

\newcommand{\CR}{\mathcal{C}_{\mathbb{R}}}
\newcommand{\CQ}{\mathcal{C}_{\mathbb{Q}}}
\newcommand{\CZ}{\mathcal{C}_{\mathbb{Z}}}
\newcommand{\Ch}{\mathcal{C}^{\mathrm{h}}}

\newcommand{\bCh}{\overline{\mathcal{C}}^{\mathrm{h}}}

\newcommand{\xl}{x}

\newcommand{\AAA}{\mathcal{A}}

\newcommand{\WWW}{\mathrm{W}}
\newcommand{\DDD}{D}

\newcommand{\CCC}{\mathrm{C}}
\newcommand{\BBB}{\mathrm{B}}
\newcommand{\HHH}{\mathrm{H}}
\newcommand{\QQQ}{\mathrm{Q}}
\newcommand{\ZZZ}{\mathrm{Z}}

\newcommand{\dG}{d_{\scl_G}}
\newcommand{\dGN}{d_{\scl_{G,N}}}
\newcommand{\dGL}{d_{\scl_{G,L}}}

\newcommand{\Gg}{G}
\newcommand{\Ng}{N}

\newcommand{\bD}{\check{D}}

\newcommand{\ta}{\tilde{a}}
\newcommand{\tb}{\tilde{b}}
\newcommand{\tc}{\tilde{c}}

\newcommand{\tg}{\tilde{g}}
\newcommand{\tx}{\tilde{x}}
\newcommand{\tz}{\tilde{z}}

\newcommand{\bT}{\overline{T}}

\newcommand{\cx}{\check{x}}
\newcommand{\cy}{\check{y}}

\newcommand{\WW}{\mathcal{W}}
\newcommand{\EE}{\mathcal{E}}
\newcommand{\FF}{\mathcal{F}}

\newcommand{\KKer}{\mathbf{Ker}}
\newcommand{\bscl}{\overline{\mathrm{scl}}}

\newcommand{\tsigma}{\tilde{\sigma}}

\newcommand{\bsim}{\pmb{\thicksim}\!\!}
\newcommand{\bbsim}{\pmb{\thicksim}}

\newcommand{\nlim}[1]{{#1}\text{-}\!\!\!\operatorname*{lim}}

\newcommand{\emb}{\Theta}

\DeclareMathOperator{\Aut}{Aut}
\DeclareMathOperator{\Ker}{Ker}
\DeclareMathOperator{\Coker}{Coker}
\DeclareMathOperator{\IA}{IA}
\DeclareMathOperator{\Imm}{Im}
\DeclareMathOperator{\Inn}{Inn}
\DeclareMathOperator{\Mod}{Mod}
\DeclareMathOperator{\scl}{scl}
\DeclareMathOperator{\cl}{cl}
\DeclareMathOperator{\sspan}{span}
\DeclareMathOperator{\Out}{Out}

\makeatletter
\newsavebox{\@brx}
\newcommand{\llangle}[1][]{\savebox{\@brx}{\(\m@th{#1\langle}\)}%
  \mathopen{\copy\@brx\kern-0.5\wd\@brx\usebox{\@brx}}}
\newcommand{\rrangle}[1][]{\savebox{\@brx}{\(\m@th{#1\rangle}\)}%
  \mathclose{\copy\@brx\kern-0.5\wd\@brx\usebox{\@brx}}}
\makeatother


\title[Coarse geometry of mixed scl II]{Coarse geometry of stable mixed commutator length II: coarse kernels for nilpotent quotients}

\author[M. Kawasaki]{Morimichi Kawasaki}
\address[Morimichi Kawasaki]{Department of Mathematics, Faculty of Science, Hokkaido University, North 10, West 8, Kita-ku, Sapporo, Hokkaido 060-0810, Japan}
\email{kawasaki@math.sci.hokudai.ac.jp}

\author[M.~Kimura]{Mitsuaki Kimura}
\address[Mitsuaki Kimura]{Department of Mathematics, Osaka Dental University, 8-1 Kuzuha Hanazono-cho, Hirakata, Osaka 573-1121, Japan}
\email{kimura-m@cc.osaka-dent.ac.jp}

\author[S. Maruyama]{Shuhei Maruyama}
\address[Shuhei Maruyama]{School of Mathematics and Physics, College of Science and Engineering, Kanazawa University, Kakuma-machi, Kanazawa, Ishikawa, 920-1192, Japan}
\email{smaruyama@se.kanazawa-u.ac.jp}

\author[T. Matsushita]{Takahiro Matsushita}
\address[Takahiro Matsushita]{Department of Mathematical Sciences, Faculty of Science, Shinshu University, Matsumoto, Nagano, 390-8621, Japan}
\email{matsushita@shinshu-u.ac.jp}

\author[M. Mimura]{Masato Mimura}
\address[Masato Mimura]{Mathematical Institute, Tohoku University, 6-3, Aramaki Aza-Aoba, Aoba-ku, Sendai 980-8578, Japan}
\email{m.masato.mimura.m@tohoku.ac.jp}

\begin{document}

\begin{abstract}
This paper studies the global difference between the stable commutator length $\scl_G$ and the stable mixed commutator length $\scl_{G,N}$ from the perspective of coarse geometry, as well as the space $\WWW(G,N)$ of non-extendable quasimorphisms modulo invariant homomorphisms. Our main theorem concerns a finitely generated group $G$ and its normal subgroup $N$ such that the quotient group $ \Gamma = G/N$ is nilpotent. 
In this setting, we prove that the natural inclusion map  of coarse groups from $[G,N]$ to $[G,G]$,  equipped with the metrics induced by $\scl_{G,N}$ and $\scl_G$, respectively, is coarsely surjective.  Moreover, this map has  a splitting coarse kernel that is coarsely isomorphic to $(\mathbb{Z}^n, \| \cdot \|_1)$, where $n$ is the dimension of ${\rm W}(G,N)$. As a consequence, $[G,N]$ is isomorphic to $\mathbb{Z}^n \times [G,G]$ as a coarse group.
\end{abstract}

\maketitle

\section{Introduction} \label{section introduction}

\subsection{Stable mixed commutator length}
The stable commutator length $\scl_G$ of a group $G$ is a fundamental invariant in geometric group theory and low-dimensional topology, bridging the algebraic properties of groups with the geometry of surfaces \cite{Calegari}. For a normal subgroup $N$ of $G$, the stable mixed commutator length $\scl_{G,N}$  was introduced in \cite{KK} as a generalization of $\scl_G$ in the sense that $\scl_{G,G}=\scl_G$.  The main theme of this paper and of Part~I \cite{fa} is to study the global difference between these two invariants. In Part~I, we characterized when $\scl_G$ and $\scl_{G,N}$ are bi-Lipschitz equivalent; in the present paper, Part~II, we describe their difference more precisely. Part~II can be read independently of Part~I.

The goal of this paper is to show that when $G$ is finitely generated and $\Gamma=G/N$ is nilpotent, the global difference between $\scl_{G,N}$ and $\scl_{G}$ can be described by the dimension of the space $\WWW(G,N)$. Here, $\WWW(G,N)$ is the space of non-extendable quasimorphisms modulo invariant homomorphisms, which is closely related to the second group cohomology and comparison maps. {The precise definition of $\WWW(G,N)$ will be given later (see equation \eqref{equation W} below).}

To explain the background of this work, we begin by reviewing the notion of invariant quasimorphisms. A \emph{quasimorphism} on a group $G$ is a function $\phi \colon G \to \RR$ such that
\[
D(\phi) = \sup_{g,h \in G} |\phi(gh) - \phi(g) - \phi(h)| < \infty.
\]
The quantity $D(\phi)$ is called the \emph{defect} of $\phi$. In addition, if $\phi(g^n) = n \cdot \phi(g)$ holds for all $g \in G$ and $n \in \ZZ$, then $\phi$ is called a \emph{homogeneous quasimorphism}.

The setting we consider in this paper is as follows. Let $N$ be a normal subgroup of $G$. A homogeneous quasimorphism $\mu$ on $N$ is said to be \emph{$G$-invariant} if $\mu(gxg^{-1}) = \mu(x)$ holds for all $g \in G$ and $x \in N$.  
Invariant quasimorphisms naturally arise in the context of symplectic geometry (see \cite{EP03, Py06, Shelukhin, FOOO},
and see \cite{E14}, \cite[Section 2]{survey} for surveys), and their invariance is applied in some papers (see \cite{EP09, commuting, Reznikov, bilinear} for example).

In the ordinary theory of quasimorphisms, it is well known that quasimorphisms are closely related to stable commutator length. Here, the \emph{commutator length} $\cl_G(x)$ of an element $x \in [G,G]$ is the minimum number of simple commutators needed to express $x$ as their product. The \emph{stable commutator length}, or \emph{scl}, is its stabilization, namely,
\[
\scl_G(x) = \lim_{n \to \infty} \frac{\cl_G(x^n)}{n}.
\]
Let $\QQQ(G)$ denote the real vector space of all homogeneous quasimorphisms on $G$. 
The celebrated Bavard duality theorem \cite{Bavard} states the following:
\begin{theorem}[\cite{Bavard}]
For every $x \in [G,G]$, the following equality holds: 
\[
\scl_G(x) = \sup_{\mu \in \QQQ(G)\setminus \HHH^1(G)} \frac{|\mu(x)|}{2D(\mu)}.
\]
\end{theorem}
The Bavard duality theorem can be generalized to the $G$-invariant setting. In this context, the counterpart of $\scl_G$ is the stabilization of the mixed commutator length $\cl_{G,N}$ defined below.

Let $N$ be a normal subgroup of $G$. For $g \in G$ and $x \in N$, an element of the form $[g,x] = gxg^{-1}x^{-1}$ is called a \emph{simple $(G,N)$-commutator}, or a \emph{simple mixed commutator}. We denote by $[G,N]$ the subgroup of $G$ generated by all simple $(G,N)$-commutators.
For $z \in [G,N]$, the \emph{$(G,N)$-commutator length} $\cl_{G,N}(z)$ of $z$ is the smallest number of simple $(G,N)$-commutators required to express $z$ as their product.
As in the ordinary case, we define $\scl_{G,N}$ as the stabilization of $\cl_{G,N}$.

Let $\QQQ(N)^G$ denote the vector space of $G$-invariant homogeneous quasimorphisms on $N$, and let $\HHH^1(N)^G$ denote the vector space of all $G$-invariant homomorphisms from $N$ to $\RR$. The following theorem generalizes Bavard's duality to {stable} mixed commutator length.

\begin{theorem}[{\cite{mixed_Bavard}}] \label{theorem mixed Bavard}
Let $G$ be a group and $N$ a normal subgroup of $G$. For every $x \in [G,N]$, the following equality holds:
\[
\scl_{G,N}(x) = \sup_{\phi \in \QQQ(N)^G\setminus \HHH^1(N)^G} \frac{|\phi(x)|}{2D(\phi)}.
\]
\end{theorem}

Thus, $\scl_{G,N}$ and $\QQQ(N)^G$ exhibit a correspondence parallel to that between $\scl_G$ and $\QQQ(G)$.

\subsection{Bi-Lipschitz equivalence and the main theorem of Part~I}

In what follows, let $G$ be a group and $N$ a normal subgroup of $G$. The main subject of this paper is the difference between the stable commutator length $\scl_G$ and the {stable mixed} commutator length $\scl_{G,N}$. More specifically, we investigate how much these two quantities can differ.

At first glance, one might expect that $\scl_G$ and $\scl_{G,N}$ typically differ. However, as the following theorem shows, they may coincide under appropriate conditions, indicating that the difference between them is a subtle phenomenon.

\begin{theorem}[{\cite[Theorems~1.10 and 2.1]{space}}]
Let $G$ be a group and $N$ a normal subgroup of $G$ such that $\HHH^2(G; \RR) = 0$ and the quotient group $ \Gamma =G/N$ is solvable. Then, for every $x \in [G,N]$, we have $\scl_{G,N}(x) = \scl_G(x)$.
\end{theorem}

This phenomenon is closely related to the space $\WWW(G,N)$ of non-extendable invariant quasimorphisms, 
whose definition is motivated by the following observation.

In general, it is known that homogeneous quasimorphisms on $G$ are conjugation-invariant \cite{Calegari}. Therefore, if a homogeneous quasimorphism $\mu$ on $N$ extends to $G$ as a homogeneous quasimorphism, then $\mu$ is necessarily $G$-invariant. On the other hand, the fact that a homogeneous quasimorphism $\mu$ on $N$ is $G$-invariant does not guarantee that it extends to $G$.

Let $i^* \colon \QQQ(G) \to \QQQ(N)^G$ be the map induced by the inclusion $i \colon N \hookrightarrow G$. Then $i^* \QQQ(G)$ is the subspace of $\QQQ(N)^G$ consisting of extendable homogeneous quasimorphisms on $N$. We define the space $\WWW(G,N)$ by
\begin{align} \label{equation W}
\WWW(G,N) = \QQQ(N)^G / (\HHH^1(N)^G + i^* \QQQ(G)).
\end{align}

We note that this space $\WWW(G,N)$ is closely related to the second ordinary cohomology of groups and the comparison maps (see Theorem~\ref{theorem VW isomorphism}), and that in several cases the dimension of $\WWW(G,N)$ is computable.  Notably,  if $G$ is Gromov-hyperbolic and $\Gamma = G/N$ is amenable, then  we  showed in \cite{space} that the dimension of $\WWW(G,N)$ coincides with the rank of $p^* \colon \HHH^2(\Gamma) \to \HHH^2(G)$ induced by the quotient map $p \colon G \to \Gamma$. In \cite{space}, we  also determined  $\dim_\RR \WWW(G,N)$ in several cases when $G$ is a surface group, a free-by-cyclic group, and the fundamental group of a mapping torus of a closed surface.

The vanishing of $\WWW(G,N)$ is closely related to the bi-Lipschitz equivalence between $\scl_G$ and $\scl_{G,N}$. Indeed, it was shown in \cite{space} that $\WWW(G,N)=0$ implies this equivalence, and in our Part~I paper \cite{fa} we established the converse, thereby obtaining the following characterization. Note that $\scl_G \le \scl_{G,N}$ always holds.

\begin{theorem}[{\cite[Theorem~2.1]{space} and \cite[Theorem~A]{fa}}] \label{theorem bi-Lipschitz}
Let $G$ be a group and $N$ a normal subgroup of $G$. Then, the following two conditions are equivalent:
\begin{enumerate}[$(1)$]
\item $\scl_G$ and $\scl_{G,N}$ are bi-Lipschitz equivalent, i.e., there exists a constant $C \ge 1$ such that for every $x \in [G,N]$,
\[
\scl_{G,N}(x) \le C \cdot \scl_G(x).
\]

\item $\WWW(G,N) = 0$.
\end{enumerate}
\end{theorem}

Theorem~\ref{theorem bi-Lipschitz} characterizes when $\scl_G$ and $\scl_{G,N}$ are bi-Lipschitz equivalent, but does not describe \emph{how} they differ when $\WWW(G,N)\ne0$. Such a description was obtained in our Part~I paper \cite{fa} in the language of coarse groups, but only at the level of rational chains, that is, of certain rational linear combinations of elements of $N$. The aim of the present paper is to carry out the description at the level of group elements, for a finitely generated group $G$ with $\Gamma=G/N$ nilpotent.

\subsection{Coarse groups and main result}

To describe the difference between $\scl_G$ and $\scl_{G,N}$, we use the theory of coarse groups developed by Leitner and Vigolo \cite{LV}.

Roughly speaking, a coarse group is a group object in the category of coarse spaces. A coarse space is an abstraction of a metric space that formalizes its large-scale structure by focusing on global connectivity rather than local distances \cite{Roe}. A typical example of a coarse group is provided by a group equipped with a bi-invariant metric, viewed as a coarse space.

Using the stable mixed commutator length $\scl_{G,N}$, we can construct a coarse group as follows. 
Define $\dGN \colon [G,N]\times[G,N] \to \mathbb{R}_{\geq 0}$ by
\[ \dGN(x,y) = \scl_{G,N}(x^{-1}y). \]
Although $\dGN$ is not a metric, it admits a bounded modification that is a genuine bi-invariant metric. Since this modification leaves the coarse structure unchanged, we may work with $([G,N], \dGN)$ as a coarse group (see Remark~\ref{remark d_scl d_scl+}).

We write $d_{\scl_{G,G}}$ simply as $\dG$. 
Then the natural inclusion map
\[
\iota^G_N \colon ([G,N], \dGN) \to ([G,G], \dG), \quad x \mapsto x
\]
is a coarse group homomorphism. In this way, the problem of comparing $\scl_G$ and $\scl_{G,N}$ is reformulated as that of comparing $\dG$ and $\dGN$.

The coarse kernel of $\iota_N^G$ measures its failure to be a
coarse embedding. We study this coarse homomorphism under the
assumptions that $G$ is finitely generated and that
$\Gamma=G/N$ is nilpotent. Under these assumptions,
$\WWW(G,N)$ is finite-dimensional
(Proposition~\ref{proposition fg nilpotent}). Our main theorem
shows that $\iota_N^G$ is coarsely surjective and that the
dimension of $\WWW(G,N)$ determines its coarse kernel.

\begin{mthm}[main theorem] \label{theorem A}
Let $N$ be a normal subgroup of a finitely generated group $G$. Assume that $\Gamma = G/N$ is nilpotent. Set $n = \dim_\RR \WWW(G,N)$. Then:
\begin{enumerate}[$(1)$]
\item $\iota^G_N \colon ([G,N], \dGN) \to ([G,G], \dG)$ is coarsely surjective.

\item There exists a splitting coarse kernel of $\iota^G_N$ (see Definition~\ref{definition splitting}) that is coarsely isomorphic to $(\ZZ^n, \| \cdot \|_1)$.
\end{enumerate}
In particular, if $r \colon ([G,N], \dGN) \to (\ZZ^n, \| \cdot \|_1)$ is the associated retraction, then the map
\[
(r, \iota^G_N)\colon ([G,N], \dGN) \xrightarrow{\cong} (\ZZ^n, \| \cdot \|_1) \times ([G,G], \dG), \quad x \mapsto (r(x), x)
\]
is an isomorphism of coarse groups.
\end{mthm}

\begin{remark}
Theorem~\ref{theorem A} reveals not merely the failure of bi-Lipschitz equivalence between $\scl_G$ and $\scl_{G,N}$ but a finer global geometric structure of the metric space $([G,N], \dGN)$. For instance, let $A \subset [G,N]$ be a subset on which $\scl_G$ is bounded. Then $(A, \dGN)$ is  coarsely contained in the coarse kernel and hence coarsely embeds into $\ZZ^n$. In particular, the asymptotic dimension of $(A, \dGN)$ is at most $n$. Moreover, there exist subsets $A$ for which this asymptotic dimension is exactly $n$.
\end{remark}

We comment on the assumption that $\Gamma =G/N$ is nilpotent. In Bridson's universe of groups \cite{Bridson}, two contrasting classes of groups have been extensively studied: amenable groups and non-positively curved groups. 
On an amenable group, every homogeneous quasimorphism is a homomorphism
(see \cite[Proposition~2.65]{Calegari}). On a non-elementary hyperbolic group, by contrast,
the space of homogeneous quasimorphisms modulo homomorphisms is
infinite-dimensional \cite{EpsteinFujiwara}.
When the quotient $\Gamma=G/N$ is hyperbolic, Tao \cite{Tao} proved
that every $G$-invariant homogeneous quasimorphism on $N$ extends
to $G$. This implies 
 that $\WWW(G,N) = 0$, and hence $\scl_G$ and $\scl_{G,N}$ are bi-Lipschitz equivalent by Theorem \ref{theorem bi-Lipschitz}. This motivates focusing on the amenable side, where nilpotent groups form a central and fundamental class.

Note that when $\Gamma$ is hyperbolic, $\iota_N^G$ is not
necessarily coarsely surjective. Indeed, if $G$ is a non-elementary
hyperbolic group and $N$ is trivial, then $([G,N],\dGN)$ consists
of a single point, whereas $([G,G],\dG)$ is unbounded by Bavard duality,
since $\QQQ(G)\neq\HHH^1(G)$ as mentioned above.

We conclude this subsection by describing the results of \cite{fa} in more detail. There we studied $\scl_{G,N}$ by functional-analytic methods: we determined when the rational-chain analogue of $\iota^G_N$ admits a coarse kernel, and identified it with the kernel of the induced bounded operator between the completed chain spaces \cite[Theorem~B]{fa}; no assumption on $G$ or $N$ was needed. At the level of group elements, only an approximation with scaling-up is available \cite[Lemma~5.1]{fa}, which sufficed for the aim of Part~I: the implication ``not $(2)$ $\Rightarrow$ not $(1)$'' in Theorem~\ref{theorem bi-Lipschitz}.

\subsection{Proof strategy and the core criterion}
\label{subsec=outline}

The proof of Theorem~\ref{theorem A} relies on the construction of a splitting coarse kernel for the natural map associated with two consecutive terms of the lower central series:
\[
\iota^G_{N,L}\colon ([G,N],\dGN)\longrightarrow([G,L],\dGL),
\qquad N=\gamma_{q+1}(G),\quad L=\gamma_q(G),
\]
where $q\geqslant1$ and $\gamma_1(G)=G$. The \emph{core criterion}, which is formulated in terms of the \emph{core} of an invariant quasimorphism, supplies the algebraic input for this construction. The core is defined below, and the basic case of the criterion is stated as Theorem~\ref{theorem B}. This construction yields Theorem~\ref{theorem C}, the one-step statement. Composing the one-step results for successive values of $q$ gives Theorem~\ref{theorem A} when $N$ is a term of the lower central series; for a general $N$, one chooses a term of the lower central series contained in $N$ and reduces to the previous case.

To formulate the one-step statement (Theorem~\ref{theorem C}), we use the following generalization of $\WWW(G,N)$.
\begin{definition}\label{definition W}
Let $G$ be a group, and let $N$ and $L$ be normal subgroups of $G$ such that $N\leqslant L$. Let $i_{N,L}\colon N\hookrightarrow L$ denote the inclusion. We define the space $\WW(G,L,N)$ by
\[
\WW(G,L,N)=\QQQ(N)^G/\bigl(\HHH^1(N)^G+i_{N,L}^*\QQQ(L)^G\bigr).
\]
\end{definition}
In particular, $\WW(G,G,N)=\WWW(G,N)$. The stable mixed commutator length $\scl_{G,N}$ induces a norm, again denoted by $\scl_{G,N}$, on a quotient $\Ch(G,N)$ of the space of real chains on $N$; let $(\bCh(G,N),\bscl_{G,N})$ denote the completion of $(\Ch(G,N),\scl_{G,N})$, and let $\emb\colon[G,N]\to\bCh(G,N)$ be the natural map $x\mapsto x$. The map $\emb$ is a coarse embedding and a coarse homomorphism (Subsection~\ref{subsec=refinedgmBavard}). The inclusion $N\hookrightarrow L$ induces a bounded linear map
\[
\bT_{N,L}\colon\bCh(G,N)\longrightarrow\bCh(G,L).
\]
When $\WW(G,L,N)$ is finite-dimensional, Theorem~\ref{theorem W duality} gives
\[
\Ker(\bT_{N,L})^*\cong\WW(G,L,N).
\]
Thus, in this case, $\Ker(\bT_{N,L})$ is finite-dimensional with dual space $\WW(G,L,N)$, and it serves as the linear model for the coarse kernel we seek. For $N=\gamma_{q+1}(G)$ and $L=\gamma_q(G)$, the required finite-dimensionality holds when $G$ is finitely generated (Corollary~\ref{corollary fg nilpotent W}~\textup{(2)}).

To find suitable elements of $\Ker(\bT_{N,L})$, take a surjective group homomorphism $\pi\colon F\to G$ from a free group $F$, with kernel $R$, and set $M=\gamma_{q+1}(F)$ and $K=\gamma_q(F)$. Every $\mu\in\QQQ(N)^G$ admits a decomposition
\[
\pi^*\mu=h+\nu|_M,
\qquad h\in\HHH^1(M)^F,\quad\nu\in\QQQ(K)^F.
\]
The homomorphism $h$ is independent of the decomposition, although it depends on the choice of $\pi$; the two conditions in Theorem~\ref{theorem B} below do not depend on this choice. We call the restriction $h|_{R\cap M}$ the \emph{core} of $\mu$. In the basic case $q=1$, where $L=G$ and $N=[G,G]$, we have the following criterion, which is an algebraic result of independent interest.

\begin{mthm}[core criterion] \label{theorem B}
Let $G$ be a group and $\pi\colon F\to G$ a surjective group homomorphism from a free group $F$, with kernel $R$. Let $\mu\in\QQQ([G,G])^G$, $h\in\HHH^1([F,F])^F$, and $\nu\in\QQQ(F)$ satisfy
\[
\pi^*\mu=h+\nu|_{[F,F]}.
\]
Then the following hold:
\begin{enumerate}[$(1)$]
\item The quasimorphism $\mu$ extends to $G$ if and only if $h=0$.
\item The class of $\mu$ in $\WWW(G,[G,G])$ vanishes if and only if $h|_{R\cap[F,F]}=0$.
\end{enumerate}
\end{mthm}

We now return to the general case. The corresponding criteria for consecutive terms of the lower central series are proved in Proposition~\ref{proposition tsigma_1} (extendability to $L$) and Theorem~\ref{theorem core kernel} (vanishing in $\WW(G,L,N)$). By Theorem~\ref{theorem core kernel}, if $\mu$ represents a nonzero class in $\WW(G,L,N)$, there exists $\tz\in R\cap M$ with $h(\tz)\ne0$. Since $M=[F,K]$, write $\tz=\prod_{j=1}^r[\tg_j,\tx_j]$ with $\tg_j\in F$ and $\tx_j\in K$. Set $g_j=\pi(\tg_j)\in G$ and $x_j=\pi(\tx_j)\in L$. The relation $\prod_{j=1}^r[g_j,x_j]=e_G$ gives a norm-limit
\[
\alpha=\nlim{\bscl_{G,N}}_{n\to\infty}
\sum_{j=1}^r\frac{[g_j,x_j^n]}{n}
\in\Ker(\bT_{N,L}),
\qquad \mu(\alpha)=h(\tz)\ne0.
\]
Thus, the relation produces an element $\alpha$ of $\Ker(\bT_{N,L})$, and the core detects the nontriviality of $\alpha$ (Proposition~\ref{proposition kernel element} and Theorem~\ref{theorem overflow}).

Suppose now that $G$ is finitely generated, and set $l=\dim_\RR\WW(G,L,N)$. Combining this construction with the duality $\Ker(\bT_{N,L})^*\cong\WW(G,L,N)$, we obtain a basis $\alpha_1,\ldots,\alpha_l$ of $\Ker(\bT_{N,L})$ consisting of such norm-limits (Corollary~\ref{corollary basis 1}). For each $\alpha_i$, the corresponding products $\alpha_i^{(m)}\in[G,N]$ of the commutators $[g_{i,j},x_{i,j}^m]$ realize $m\alpha_i$ with an error independent of $m$ (Lemma~\ref{lemma lifting}). This is the key step in passing from chains to group elements: the approximation in \cite[Lemma~5.1]{fa} realizes only an integer multiple $k\alpha$ of a chain $\alpha$ by an element of $[G,N]$, with an error proportional to $k$, and one cannot in general divide by $k$ within $[G,N]$. Hence $(m_1,\ldots,m_l)\mapsto\alpha_1^{(m_1)}\cdots\alpha_l^{(m_l)}$ is a coarse homomorphism from $(\ZZ^l,\|\cdot\|_1)$ to $([G,N],\dGN)$ realizing $\Ker(\bT_{N,L})$ (Proposition~\ref{proposition group lifting}), and quasimorphisms representing the dual basis provide a coarse retraction (Proposition~\ref{proposition splitting}). Together with the coarse surjectivity of $\iota^G_{N,L}$ (Corollary~\ref{corollary xi equivalence}), this proves the following one-step statement.

\begin{mthm}[one-step case] \label{theorem C}
Let $G$ be a finitely generated group and $q\geqslant1$ an integer. Set
\[
N=\gamma_{q+1}(G),\qquad L=\gamma_q(G),\qquad
l=\dim_\RR\WW(G,L,N).
\]
Then the natural map
\[
\iota^G_{N,L}\colon([G,N],\dGN)\longrightarrow([G,L],\dGL)
\]
is coarsely surjective and has a splitting coarse kernel coarsely isomorphic to $(\ZZ^l,\|\cdot\|_1)$.
\end{mthm}

For $q=1$ we have $L=G$ and $N=[G,G]$, so that $\iota^G_{N,L}=\iota^G_{[G,G]}$ and $l=\dim_\RR\WWW(G,[G,G])$; in this case Theorem~\ref{theorem C} gives assertions $(1)$ and $(2)$ of Theorem~\ref{theorem A} for $N=[G,G]$.

To deduce Theorem~\ref{theorem A}, choose $k$ with $\gamma_k(G)\leqslant N$. Composing the one-step maps $\iota^G_{\gamma_{i+1}(G),\gamma_i(G)}$ $(1\leqslant i<k)$, each of which is coarsely surjective (Theorem~\ref{theorem C}), we see that $\iota^G_{\gamma_k(G)}$, and hence $\iota^G_N$, is coarsely surjective, and that $\iota^G_{\gamma_k(G)}$ has a splitting coarse kernel coarsely isomorphic to $(\ZZ^{m},\|\cdot\|_1)$, where $m=\dim_\RR\WWW(G,\gamma_k(G))$ (Theorem~\ref{theorem direct sum}). We then extract the part of this coarse kernel corresponding to $\Ker(\bT_{N,G})$ (Lemmas~\ref{lemma surjective} and~\ref{lemma factorization}), which gives a splitting coarse kernel of $\iota^G_N$ coarsely isomorphic to $(\ZZ^n,\|\cdot\|_1)$. See Section~\ref{section main theorem} for details.

The construction just described is not explicit: in general, there is no systematic way to choose relations $\tz \in R \cap M$ whose cores detect the nonzero classes in $\WW(G,L,N)$, and hence the coarse kernel is not obtained as an explicit subset of $[G,N]$. In Section~\ref{section explicit}, we nevertheless determine the coarse kernel explicitly when suitable relations are available, namely for surface groups, free-by-cyclic groups, and surface mapping-torus groups. For instance, if $G = \langle a_1, \ldots, a_l, b_1, \ldots, b_l \mid [a_1,b_1] \cdots [a_l,b_l] \rangle$ with $l \geq 2$ and $N = [G,G]$, then $\dim_\RR \WWW(G,N) = 1$ and the set $\{ [a_1, b_1^m] \cdots [a_l, b_l^m] \mid m \in \ZZ \}$ realizes the coarse kernel of $\iota^G_N$. The explicit coarse kernels also give concrete criteria for the vanishing of a class in $\WWW(G,N)$: a $G$-invariant homogeneous quasimorphism on $N$ represents the zero class precisely when the restriction of the quasimorphism to the coarse kernel is bounded (Proposition~\ref{proposition W vanishing}).

\subsection{Organization of the paper}

The rest of this paper is organized as follows: In Section~\ref{section preliminaries}, we review the preliminary concepts of coarse groups, invariant quasimorphisms, and stable commutator lengths. In Section~\ref{section functional analysis}, we show that the continuous dual of $\Ker(\bT_{N,L})$ is isomorphic to $\WW(G,L,N)$, under mild conditions. In Section~\ref{section kernel element}, we construct certain elements of $\Ker(\bT_{N,L})$ under certain conditions on $G$, $L$ and $N$. In Section~\ref{section core}, we introduce the concept of the core and discuss the coarse kernel of $\iota^G_{N,L} \colon ([G,N], \dGN) \to ([G,L], \dGL)$ in the case where $L =  \gamma_{q} (G)$ and $N =  \gamma_{q+1} (G)$. In Section~\ref{section main theorem}, we complete the proof of Theorem~\ref{theorem A}. In Section~\ref{section explicit}, we provide explicit representations of the coarse kernel of $\iota^G_N$ for several examples.

\subsection{Notation and conventions} \label{subsection notation}

In the rest of this paper, $G$ denotes a group, and $L$ and $N$ denote normal subgroups of $G$ such that $N \leqslant L$. Here, the symbol $\leqslant$ denotes that $N$ is a subgroup of $L$. Let $i_{N, L} \colon N \hookrightarrow L$ denote the inclusion. Let $\mathbb{N}$ denote the set of positive integers. The symbol $\gamma_q(G)$ denotes the $q$-th term of the lower central series of $G$, defined inductively by $ \gamma_1(G)  = G$ and $\gamma_{q+1}(G) = [G, \gamma_q(G)]$. If the coefficient group is omitted in cohomology and bounded cohomology, it is assumed to be $\RR$. For example, $\HHH^n(G)$ means $\HHH^n(G; \RR)$, and $\HHH^n_b(G)$ means $\HHH^n_b(G; \RR)$. For $a,b \in \RR$ and $D\geq 0$, we write $a \sim_D b$ to mean $|a-b| \leq D$.

\section{Preliminaries} \label{section preliminaries}

In this section, we review several definitions and facts we need in this paper, following \cite{Calegari, Frigerio}.

\subsection{(Bounded) cohomology of groups}

Let $G$ be a group and let $n \in \mathbb{Z}$. For $n \ge 0$, let $\CCC^n(G)$ be the $\RR$-module of all $\RR$-valued functions on $G^n$, and set $\CCC^n(G) = 0$ for $n < 0$. 
Define the coboundary map $\delta \colon \CCC^n(G) \to \CCC^{n+1}(G)$ by
\begin{multline*}
\delta c(g_1, \cdots, g_{n+1}) 
= c(g_2, \cdots, g_{n+1})
+ \sum_{i=1}^n (-1)^i c(g_1, \cdots, g_i g_{i+1}, \cdots, g_{n+1}) \\
+ (-1)^{n+1} c(g_1, \cdots, g_n)
\end{multline*}
for $n \ge 1$, and set $\delta = 0$ for $n \le 0$. The \emph{$n$-th group cohomology $\HHH^n(G)$ of $G$} is the cohomology of the complex $(\CCC^{\bullet} (G), \delta)$ in degree $n$. In particular, $\HHH^1(G)$ is canonically identified with $\mathrm{Hom}(G, \RR)$.

Let $\CCC^n_b(G)$ be the space of bounded $\RR$-valued functions on $G^n$ for $n \ge 0$, and set $\CCC^n_b(G) = 0$ for $n < 0$. Then $\CCC^\bullet_b(G)$ is a subcomplex of $\CCC^\bullet(G)$, and its $n$-th cohomology $\HHH^n_b(G)$ is called the \emph{$n$-th bounded cohomology} of $G$. The inclusion $\CCC^\bullet_b(G) \hookrightarrow \CCC^\bullet(G)$ induces the \emph{comparison map} $c_G^\bullet \colon \HHH^\bullet_b(G) \to \HHH^\bullet(G)$.

\subsection{Invariant quasimorphisms}\label{subsec=qm}

A \emph{quasimorphism on $G$} is a function $\phi \colon G \to \RR$ such that
\[
D(\phi) = \sup_{g,h \in G} |\phi(gh) - \phi(g) - \phi(h)|
\]
is finite. The value $D(\phi)$ is called the \emph{defect of $\phi$}. A quasimorphism $\phi$ is said to be \emph{homogeneous} if $\phi(g^n) = n \phi(g)$ for every $g \in G$ and $n \in \ZZ$. Let $\QQQ(G)$ denote the $\RR$-linear space consisting of homogeneous quasimorphisms on $G$.

A homogeneous quasimorphism $\mu \colon N \to \RR$ is called \emph{$G$-invariant} if $\mu(gxg^{-1}) = \mu(x)$ for every $g \in G$ and $x \in N$. Let $\QQQ(N)^G$ denote the $\RR$-linear space consisting of $G$-invariant homogeneous quasimorphisms on $N$.

A homogeneous quasimorphism $\mu$ on $N$ is said to be \emph{extendable to $G$} if there exists $\phi \in \QQQ(G)$ such that $\phi|_N = \mu$. 
Since every homogeneous quasimorphism on $G$ is $G$-invariant (see \cite{Calegari} for example), every homogeneous quasimorphism on $N$ that is extendable to $G$ is $G$-invariant.

Let $i_{N,L}^* \colon \QQQ(L)^G \to \QQQ(N)^G$ be the map defined by restriction,
\[
i_{N,L}^* \colon \QQQ(L)^G \to \QQQ(N)^G,\quad \phi \mapsto \phi|_N.
\]
Note that $i_{N,G}^* \QQQ(G)^G = i_{N,G}^* \QQQ(G)$ is the $\RR$-linear space consisting of homogeneous quasimorphisms extendable to $G$. 
Recall from Definition~\ref{definition W} that
\[
\WW(G, L, N) = \QQQ(N)^G / (\HHH^1(N)^G + i_{N, L}^* \QQQ(L)^G).
\]
We note that the $\RR$-linear map $i_{N,L}^*\colon \QQQ(L)^G \to \QQQ(N)^G$ induces an $\RR$-linear map
\[
\widehat{i^*}_{N,L}\colon \QQQ(L)^G/\HHH^1(L)^G \to \QQQ(N)^G/\HHH^1(N)^G,
\]
and that $\WW(G, L, N)$ can be regarded as the cokernel of $\widehat{i^*}_{N,L}$.

We express $D(\phi)$ as $D_L(\phi)$ for $\phi\in \QQQ(L)^G$ and $D(\mu)$ as $D_N(\mu)$ for $\mu\in \QQQ(N)^G$, respectively. Then, the defect $D_L$, as a map from $\QQQ(L)^G$ to $\RR_{\geq 0}$, gives rise to a genuine norm on the $\RR$-linear space $\QQQ(L)^G/\HHH^1(L)^G$, and we write this norm as $D_L$ again. We define the defect norm $D_N$ on $\QQQ(N)^G/\HHH^1(N)^G$ similarly. Then, $\widehat{i^*}_{N,L}$ may be regarded as a map between two real normed spaces $(\QQQ(L)^G/\HHH^1(L)^G,D_L)$ and $(\QQQ(N)^G/\HHH^1(N)^G,D_N)$. In particular, if $\widehat{i^*}_{N,L}(\QQQ(L)^G/\HHH^1(L)^G)$ is norm-closed in $(\QQQ(N)^G/\HHH^1(N)^G,D_N)$, then  $\WW(G, L, N)$ is equipped with the quotient norm (we recall its definition in Subsection~\ref{subsec=fa}), which we write $\bD_N$. 

A group $\Gamma$ is said to be \emph{boundedly $3$-acyclic} if $\HHH^k_b(\Gamma) = 0$ for $k = 1, 2, 3$. Amenable groups are typical examples of boundedly $3$-acyclic groups.

\begin{theorem}[{\cite[Theorem~1.10]{space} and \cite[Corollary~7.20]{KIAS}}] \label{theorem VW isomorphism}
Let $G$ be a group, $N$ a normal subgroup of $G$, and $p \colon G \to \Gamma = G/N$ the quotient map. Suppose that $\Gamma$ is boundedly $3$-acyclic. Then there is a natural isomorphism
\[
\WWW(G,N) \cong \Imm(c_G^2 \colon \HHH^2_b(G) \to \HHH^2(G)) \cap \Imm(p^* \colon \HHH^2(\Gamma) \to \HHH^2(G)).
\]
\end{theorem}

Therefore, we obtain the following corollary.

\begin{corollary} \label{corollary decomposition}
Let $F$ be a free group and let $M$ be a normal subgroup of $F$. If $F / M$ is boundedly $3$-acyclic, then $\WWW(F, M) = 0$. In other words, for every $F$-invariant homogeneous quasimorphism $\mu$ on $M$, there exist $f \in \HHH^1(M; \RR)^F$ and $\phi \in \QQQ(F)$ such that $\mu = f + \phi|_{M}$.
\end{corollary}

Many theorems in this paper rely on the assumption that $\WW(G, L, N)$ is finite-dimensional. Here, we give criteria for $\WW(G, L, N)$ to be finite-dimensional.

The following proposition follows from Theorem~\ref{theorem VW isomorphism} and the finite-dimensionality of the real cohomology groups of finitely generated nilpotent groups.

\begin{proposition}[{see \cite[Corollary~7.18]{KIAS}}] \label{proposition fg nilpotent}
Let $G$ be a finitely generated group, $N$ a normal subgroup of $G$ such that $\Gamma = G/N$ is nilpotent. Then $\WWW(G,N)$ is finite-dimensional.
\end{proposition}

\begin{corollary} \label{corollary fg nilpotent W}
Let $G$ be a group, and let $N$ and $L$ be normal subgroups of $G$ such that $N \leqslant L$. Then the following hold.
\begin{enumerate}[$(1)$]
\item The space $\WW(G,L,N)$ is the cokernel of the natural map $\WWW(G,L) \to \WWW(G,N)$, $[\mu] \mapsto [\mu|_N]$. In particular, if $\WWW(G,N)$ is finite-dimensional, then so is $\WW(G,L,N)$.

\item If $G$ is finitely generated and $\Gamma = G/N$ is nilpotent, then $\WW(G,L,N)$ is finite-dimensional.
\end{enumerate}
\end{corollary}
\begin{proof}
Since every homogeneous quasimorphism on $G$ restricts to a $G$-invariant homogeneous quasimorphism on $L$, we have $i^* \QQQ(G) \subseteq i_{N,L}^* \QQQ(L)^G$ in $\QQQ(N)^G$. Hence \textup{(1)} follows directly from the definitions of $\WWW(G,N)$, $\WWW(G,L)$ and $\WW(G,L,N)$. Then \textup{(2)} follows from \textup{(1)} and Proposition~\ref{proposition fg nilpotent}.
\end{proof}

\begin{remark}
Suppose that $\Gamma = G/N$ is boundedly $3$-acyclic and that at least one of $\HHH^2_b(G)$, $\HHH^2(G)$, and $\HHH^2(\Gamma)$ is finite-dimensional. Then $\WWW(G,N)$ is finite-dimensional by Theorem~\ref{theorem VW isomorphism}, and hence so is $\WW(G,L,N)$ for every normal subgroup $L$ of $G$ with $N \leqslant L$ by Corollary~\ref{corollary fg nilpotent W}~\textup{(1)}.
\end{remark}

We have the following additive dimension formula. 

\begin{proposition} \label{proposition W injection}
Let $N$ and $L$ be normal subgroups of a group $G$ such that $N \leqslant L$. Assume that $G/N$ and $G/L$ are boundedly $3$-acyclic. Then the following hold:
\begin{enumerate}[$(1)$]
\item The map $\WWW(G,L) \to \WWW(G,N)$, $[\mu] \mapsto [\mu|_N]$ is injective.
\item Let $L'$ be a normal subgroup of $G$ with $L \leqslant L'$ such that $G/L'$ is boundedly $3$-acyclic. If $\WWW(G,N)$ is finite-dimensional, then
\[ \dim_\RR \WW(G,L',N) = \dim_\RR \WW(G,L',L) + \dim_\RR \WW(G,L,N).\]
In particular, taking $L' = G$, we have $\dim_\RR \WWW(G,N) = \dim_\RR \WWW(G,L) + \dim_\RR \WW(G,L,N)$.
\end{enumerate}
\end{proposition}
\begin{proof}
Let $p_L\colon G\to G/L$ and $p_N\colon G\to G/N$ be the quotient maps. By the naturality of Theorem~\ref{theorem VW isomorphism}, we have the commutative diagram
\[ \xymatrix{
\WWW(G,L) \ar[d] \ar[rr]^-{\cong} & & \Imm(c_G^2 \colon \HHH^2_b(G) \to \HHH^2(G)) \cap \Imm(p_L^* \colon \HHH^2(G/L) \to \HHH^2(G)) \ar[d] \\
\WWW(G, N) \ar[rr]^-{\cong} & & \Imm(c_G^2 \colon \HHH^2_b(G) \to \HHH^2(G)) \cap \Imm(p_N^* \colon \HHH^2(G/N) \to \HHH^2(G)).
}\]
Since $G \to G/L$ factors through $G/N$, the right vertical arrow is injective. This proves \textup{(1)}.

For \textup{(2)}, applying \textup{(1)} to the pairs $(N,L)$ and $(L,L')$, we see that the restriction maps $\WWW(G,L') \to \WWW(G,L) \to \WWW(G,N)$ are injective, and their composition is the restriction map $\WWW(G,L') \to \WWW(G,N)$. By Corollary~\ref{corollary fg nilpotent W}~\textup{(1)}, the spaces $\WW(G,L',N)$, $\WW(G,L',L)$ and $\WW(G,L,N)$ are the cokernels of $\WWW(G,L') \to \WWW(G,N)$, $\WWW(G,L') \to \WWW(G,L)$ and $\WWW(G,L) \to \WWW(G,N)$, respectively. Since $\WWW(G,N)$ is finite-dimensional, the formula follows.
\end{proof}

\subsection{Stable mixed commutator length}\label{subsec=mixedscl}

For a comprehensive survey of stable mixed commutator lengths, we refer the reader to \cite{survey, KIAS}. A \emph{simple mixed commutator}, or \emph{simple $(G,N)$-commutator}, is an element of $G$ of the form $[g,x] = gxg^{-1}x^{-1}$, where $g \in G$ and $x \in N$. Let $[G,N]$ denote the subgroup generated by the set of all simple $(G,N)$-commutators. Then $[G,N]$ is a normal subgroup of $G$ contained in $N$. The \emph{mixed commutator length}, or \emph{$(G,N)$-commutator length}, denoted by $\cl_{G,N}$, is defined to be the word length on $[G,N]$ with respect to the set of simple $(G,N)$-commutators. The \emph{stable mixed commutator length}, or \emph{stable $(G,N)$-commutator length} of $y \in [G,N]$ is defined by
\[
\scl_{G,N}(y) = \lim_{n \to \infty} \frac{\cl_{G,N}(y^n)}{n}.
\]
Here, Fekete's lemma guarantees the existence of the limit on the right-hand side.

Let $A = \ZZ$, $\QQ$ or $\RR$. Let  $\CCC_n(G;A)$  be the $A$-module freely generated by  $G^n = G\times \dots \times G$.  Define the \emph{boundary map} $\partial \colon \CCC_2(G; A)\to \CCC_1(G; A)$ to be the $A$-linear map satisfying
\[ \partial ((g_1, g_2))=g_2- g_1 g_2+g_1\]
for every $(g_1, g_2) \in G^2$. Set $\ZZZ_2(G;A) = \Ker \partial$ and $\BBB_1(G; A) = \partial \CCC_2(G; A)$. Let $\CCC_2'(G,N; A)$ be the $A$-submodule of $\CCC_2(G;A)$ generated by the set
\[ \{ (g_1, g_2) \mid \text{one of $g_1$ and $g_2$ belongs to $N$}\}.\]
Set $\ZZZ'_2(G,N; A) = \CCC_2'(G, N; A) \cap \ZZZ_2(G; A)$ and $\BBB'_1(G, N; A) = \partial \CCC_2'(G,N; A)$. We regard $\CCC_2(G; \RR)$ and $\CCC_2'(G,N; \RR)$ as normed spaces with $\ell^1$-norm. Then $\ZZZ_2'(G, N; \RR)$ is a closed subspace of $\CCC_2'(G, N ; \RR)$ and we regard $\BBB_1'(G,N; \RR)$ as the normed space $\CCC_2'(G,N ; \RR) / \ZZZ_2'(G,N ; \RR)$ with the quotient norm
  $\|\cdot\|'$. Using the natural inclusion
$\CCC_1(N;A)\subset\CCC_1(G;A)$, define
$\CC_A(G,N) = \CCC_1(N;A)\cap\BBB'_1(G,N;A)$.
We equip $\CR(G,N)$ with the restriction of $\|\cdot\|'$, which we write $\|\cdot\|'$ again.

\begin{lemma}[{\cite{KIAS}}]\label{lem=domainG}
Let $G$ be a group, $N$ a normal subgroup of $G$ and let $c \in \CCC_1(N ; \ZZ)$. Then the following are equivalent:
\begin{enumerate}
 \item[\textup{(i)}] There exist $k,l\in \ZZ_{\geq 0}$ and $x_1, \cdots, x_k, \cx_1, \cdots, \cx_l \in N$ such that
\[
c = x_1 + \cdots + x_k - \cx_1 - \cdots - \cx_l
\]
and $x_1 \cdots x_k \cx_1^{-1} \cdots \cx_l^{-1} \in [G,N]$.
 \item[\textup{(ii)}] $c$ belongs to $\BBB'_1(G, N;\ZZ)$.
\end{enumerate}
\end{lemma}

We now define $\scl_{G,N}$ on $\CZ(G,N)$.

\begin{definition}[{\cite[Definitions~5.25 and 5.29]{KIAS}}] \label{scl for integral}
Let $c\in\CZ(G,N)$ and write $c = x_1 + \cdots + x_k - \cx_1 - \cdots - \cx_l$ as in Lemma~\ref{lem=domainG} \textup{(i)}, with $x_i\ne\cx_j$ for all $i,j$.
\begin{enumerate}[$(1)$]
  \item Define
\[
\cl_{G,N}(c) = \inf \cl_{G,N}(y_1 \cdots y_k \cy_1^{-1} \cdots \cy_l^{-1}),
\]
where $y_i$ runs over the $G$-conjugacy class of $x_i$, and $\cy_j$ runs over the $G$-conjugacy class of $\cx_j$.

  \item Define
\[
\scl_{G,N}(c) = \lim_{n\to \infty} \frac{1}{n} \cl_{G,N}(x_1^n + \cdots + x_k^n  +  \cx_1^{-n}  +  \cdots  +  \cx_l^{-n}).
\]

  \item Let $c\in\CQ(G,N)$. Then there exists $m \in \NN$ such that $mc\in\CZ(G,N)$. Define
\[
\scl_{G,N}(c) = \frac{1}{m} \scl_{G,N}(mc).
\]
\end{enumerate}
\end{definition}
The existence of the limit in (2) follows from \cite[Lemma~5.19]{KIAS}, and the well-definedness of $\cl_{G,N}$ and $\scl_{G,N}$ is proved in \cite{KIAS}.

We now state the generalized mixed Bavard duality theorem (Theorem~\ref{theorem gBavard}). 
For $\mu\in\QQQ(N)^G$ and
$c=a_1x_1+\cdots+a_rx_r\in\CR(G,N)$, where
$a_1,\ldots,a_r\in\RR$ and $x_1,\ldots,x_r\in N$, set
\[
\mu(c)=a_1\mu(x_1)+\cdots+a_r\mu(x_r).
\]
Recall from Subsection~\ref{subsec=qm} that $\QQQ(N)^G/\HHH^1(N)^G$ is equipped with the defect norm $D_N$. Also, we note that for $c\in \CR(G,N)$, the value $\mu(c)$ only depends on the class $[\mu]\in \QQQ(N)^G/\HHH^1(N)^G$ by   the $\RR$-linear extension of  Lemma~\ref{lem=domainG}.

\begin{theorem}[{generalized mixed Bavard duality theorem \cite{KIAS}}] \label{theorem gBavard}
For $c\in\CQ(G,N)$, we have
\[
\scl_{G,N}(c) = \sup_{[\mu] \in \big(\QQQ(N)^G /\HHH^1(N)^G\big)\setminus \{0\}} \frac{|\mu(c)|}{2 D_N(\mu)}.
\]
\end{theorem}

We extend $\scl_{G,N}$ to $\CR(G,N)$ by the right-hand side of Theorem~\ref{theorem gBavard}. That is, for $c = a_1 x_1 + \cdots + a_k x_k \in \CR(G,N)$ with $a_1,\ldots,a_k\in\RR$ and $x_1,\ldots,x_k\in N$, we define
\[
\scl_{G,N}(c) = \sup_{[\mu] \in \big(\QQQ(N)^G /\HHH^1(N)^G\big)\setminus \{0\}} \frac{|\mu(c)|}{2 D_N(\mu)}.
\]

\subsection{A few facts on functional analysis}\label{subsec=fa}

In this subsection, we recall some basic facts on functional analysis, which will be used in Subsection~\ref{subsec=refinedgmBavard} and Section~\ref{section functional analysis}. For a real normed space $X=(X,\|\cdot\|_X)$, $X^{\ast}$ denotes its continuous dual. The pairing $\langle \cdot ,\cdot\rangle_X\colon X^{\ast}\times X\to \RR$; $(f,x)\mapsto f(x)$ is called the standard pairing for $X$. The dual norm $\|\cdot\|_{X^{\ast}}$ on $X^{\ast}$ is defined by 
\[
\|f\|_{X^{\ast}}=\sup_{x\in X\setminus \{0\}}\frac{|\langle f,x\rangle_X|}{\|x\|_X}
\]
for $f\in X^{\ast}$. For a bounded linear operator $S\colon X\to Y$ between real normed spaces, the \emph{adjoint operator} $S^{\dagger}\colon Y^{\ast}\to X^{\ast}$ is defined to be the bounded linear operator uniquely determined by the following condition: for every $x\in X$ and $f\in Y^{\ast}$,
\[
\langle S^{\dagger}f,x\rangle_X=\langle f,Sx\rangle_Y.
\]
Here, $\langle \cdot ,\cdot\rangle_X$ and $\langle \cdot ,\cdot\rangle_Y$ are the standard pairings for $X$ and $Y$, respectively.

\begin{definition}\label{def=dualitypairing}
For real normed spaces $X$ and $Y$, an $\RR$-bilinear form
\[
\langle \cdot ,\cdot\rangle\colon Y\times X\to \RR
\]
is called a \emph{duality pairing for $X$} if the map 
\[
Y \to \mathrm{Hom}_{\RR}(X,\RR);\ y\mapsto \langle y,\cdot \rangle, \quad \textrm{where}\quad \langle y,\cdot \rangle\colon X\to \RR;\ x\mapsto \langle y,x\rangle,
\]
induces a surjective map from $Y$ to $X^{\ast}$ that is isometric; namely, we have $\|\langle y,\cdot\rangle\|_{X^{\ast}}=\|y\|_Y$ for every $y \in Y$. Here, $X^{\ast}$ is naturally regarded as a subset of $\mathrm{Hom}_{\RR}(X,\RR)$. 
\end{definition}
For instance, the standard pairing $\langle \cdot,\cdot\rangle_X\colon X^{\ast}\times X\to \RR$ is a duality pairing for $X$.

Practically, it is convenient to identify the continuous dual $X^*$ with a pair of a real normed space $Y$ and a duality pairing $\langle \cdot,\cdot\rangle\colon Y\times X\to \RR$. In this formulation, $X^{\ast}$ is identified with $Y$ via the duality pairing $\langle \cdot,\cdot\rangle$: the isometric isomorphism is given by $Y\ni y\mapsto \langle y,\cdot\rangle \in X^{\ast}$. In this convention, the adjoint operator of $S\colon X \to Y$ may be formulated as follows. Let $(V,\langle \cdot,\cdot\rangle)$ and $(W,\langle \cdot,\cdot\rangle')$ be pairs that are identified with $X^{\ast}$ and $Y^{\ast}$, respectively. Then, the adjoint operator $S^{\dagger}$, with respect to the two duality pairings $\langle \cdot,\cdot\rangle$ and $\langle \cdot,\cdot\rangle'$ above, is the bounded linear operator $T\colon W\to V$ uniquely determined by the following condition: $\langle Tw,x\rangle=\langle w,Sx\rangle'$ for every $x\in X$ and $w\in W$. In other words, this $T$ is the operator that makes the following diagram commute
\begin{equation*} 
\vcenter{\xymatrix{
Y^* \ar[r]^{S^\dagger} & X^*  \\
W \ar[r]^{T} \ar[u]^{\cong} & V \ar[u]_{\cong}, 
}}
\end{equation*}
where $S^\dagger$ is defined with respect to the standard pairings; we write ``$T=S^{\dagger}$'' under the identifications $V\cong X^{\ast}$ and $W\cong Y^{\ast}$ via the duality pairings $\langle \cdot,\cdot\rangle$ and $\langle \cdot,\cdot\rangle'$.

In what follows, we recall a few definitions and results on functional analysis. For simplicity, we formulate them with respect to the standard pairing. However, similarly to the discussion above for adjoint operators, we may switch to the setting with respect to a given duality pairing.

Let $X$ be a real normed space, and $\langle \cdot,\cdot \rangle_X\colon X^{\ast}\times X\to \RR$ the standard pairing for $X$. For a subspace $U$ of $X$ and a subspace $V$ of  $X^{\ast}$, define
\begin{align*} 
U^\perp &= \{ f \in X^* \mid \langle f,x\rangle_X = 0\textrm{ for every }x\in U\}\quad \textrm{and} \\
{}^\perp V &= \{ x \in X \mid \langle f,x\rangle_X = 0\textrm{ for every }f\in V\}.
\end{align*}
Then, $U^{\perp}$ is a closed subspace of $X^{\ast}$ and ${}^{\perp}V$ is a closed subspace of $X$. For a real normed space $(X,\|\cdot\|_X)$ and a closed subspace $U$ of $X$, regard $X/U$ as the quotient space equipped with the quotient norm:
\[
\|\cdot\|_{X/U}\colon X/U\to \RR_{\geq 0};\quad \|[x]\|_{X/U}=\inf\{\|x+u\|_X\;|\;u\in U\}.
\]

\begin{lemma}[{see \cite[Theorem~4.9]{Rudin} or \cite[Lemma~3.4]{fa}}] \label{lemma duality pairing}
Let $X$ be a real normed space and $U$ a subspace. Then, the standard pairing $\langle \cdot, \cdot \rangle_X \colon X^* \times X \to \RR$ induces a well-defined map
\[ \langle \cdot , \cdot \rangle \colon (X^* / U^\perp) \times U \to \RR, \quad \langle [\phi], x \rangle = \langle \phi, x \rangle_X,\]
and this $\langle \cdot, \cdot \rangle$ is a duality pairing for $U$. In particular, $U^{\ast}$ is isometrically isomorphic to $X^* / U^\perp$. 
\end{lemma}

The following theorem is \emph{Banach's closed range theorem} for bounded linear operators.

\begin{theorem}[{see \cite{Rudin} or \cite[Theorem~3.5]{fa}}] \label{theorem CR}
Let $X$ and $Y$ be real Banach spaces and $S \colon X \to Y$ a bounded linear operator. Then  the following four conditions are all equivalent.

\begin{enumerate}[{\rm(1)}]
  \item $S(X)$ is norm-closed in $Y$.
  \item $S^{\dagger}(Y^{\ast})$ is norm-closed in $X^{\ast}$.
  \item $S(X)={}^{\perp}\Ker(S^{\dagger})$.
  \item $S^{\dagger}(Y^{\ast})=\Ker(S)^{\perp}$.
\end{enumerate}
\end{theorem}

\begin{corollary}\label{cor=Banach}
Let $X$ and $Y$ be real Banach spaces and $S \colon X \to Y$ a bounded linear operator. Assume that one of the four \textup{(}equivalent\textup{)}  conditions stated in Theorem~\textup{\ref{theorem CR}} is satisfied. Then, the following map
\[
\Coker(S^\dagger) \times \Ker(S) \to \RR;\quad ([f],x)\mapsto \langle f,x\rangle_X
\]
is well-defined, and this is a duality pairing for $\Ker(S)$. Here, $\langle \cdot,\cdot\rangle_X$ is the standard pairing for $X$. In particular, $\Ker(S)^{\ast}$ is isometrically isomorphic to $\Coker(S^\dagger)=X^{\ast}/S^{\dagger}(Y^{\ast})$.
\end{corollary}

\begin{proof}
By Theorem~\ref{theorem CR}, we have $S^{\dagger}(Y^{\ast})=\Ker(S)^{\perp}$. Apply Lemma~\ref{lemma duality pairing} for $U=\Ker(S)$. This completes our proof because $X^{\ast}/U^{\perp}=X^{\ast}/S^{\dagger}(Y^{\ast})=\Coker(S^{\dagger})$.
\end{proof}

\subsection{Refined version of the generalized mixed Bavard duality theorem}\label{subsec=refinedgmBavard}
In this subsection, we recall some definitions and a result from \cite{fa}. This result is the refined version of the generalized mixed Bavard duality theorem (Theorem~\ref{theorem duality}), which plays a key role in the functional-analytic study of $\mathrm{scl}_{G,N}$.

\begin{definition} \label{defn=Ch}
Let $\Gg$ be a group and $\Ng$ its normal subgroup. 
\begin{enumerate}[$(1)$]
  \item Define $\hh(\Gg,\Ng)$ as the linear subspace of $\CR(G,N)$ generated by 
\[
	\{\xl^k-kx \mid x\in \Ng,\, k\in \ZZ\}.
\] 
  \item Define $\Ch(G,N)$ to be
\[
\Ch(G,N)=\CR(\Gg,\Ng)/\overline{\mathrm{h}(\Gg,\Ng)}^{\|\cdot\|'}.
\]
Here $\overline{U}^{\|\cdot\|'}$ denotes the closure of $U \subset \CR(G,N)$ with respect to the norm $\|\cdot\|'$.
\end{enumerate}
Let $\| \cdot \|_\hh$ denote the quotient norm on $\Ch(G,N)$, and $\bCh(G,N)$ the completion of $(\Ch(G,N), \| \cdot \|_\hh)$. We write $\Ch(G)$ for $\Ch(G,G)$ and $\bCh(G)$ for $\bCh(G,G)$, respectively.
\end{definition}

The following theorem is proved in \cite[Theorem~4.6]{fa}; see also \cite[Subsection~5.1]{fa} for the formulation of (2). 

\begin{theorem}[{refined version of generalized mixed Bavard duality theorem}]\label{theorem duality}
Let $\Gg$ be a group and $\Ng$ a normal subgroup. Then the following hold:
\begin{enumerate}[$(1)$]
\item For every $c \in \Ch(G,N)$, we have
\[ \scl_{G,N}(c) = \frac{1}{2}\| c\|_\hh.\]
Hence $\scl_{G,N}$ has a unique extension $\bscl_{G,N}$ to $\bCh(G,N)$ as a norm.

\item The map
\[
\Big(\QQQ(\Ng)^{\Gg} / \HHH^1(\Ng)^{\Gg},2\DDD_{\Ng}\Big)\times (\bCh(G,N),\bscl_{\Gg,\Ng})\to \RR;\quad ([\mu],c) \mapsto \mu(c),
\]
is a duality pairing for $(\bCh(G,N),\bscl_{\Gg,\Ng})$. In particular, $(\bCh(G,N),\bscl_{\Gg,\Ng})^{\ast}$ is isometrically isomorphic to $(\QQQ(\Ng)^{\Gg} / \HHH^1(\Ng)^{\Gg},2\DDD_{\Ng})$.
\end{enumerate}
\end{theorem}
Theorem~\ref{theorem duality} is a refinement of Theorem~\ref{theorem gBavard} because the former immediately implies the latter with the aid of the Hahn--Banach theorem. We also note that
 $\bCh(G,N)$ is equipped with two natural norms, $\| \cdot \|_\hh$ and $\bscl_{G,N}$, which differ only by a constant factor. Unless otherwise specified, we shall take $\bscl_{G,N}$ as the norm on $\bCh(G,N)$. We consider $\Ch(G,N)$ in a similar manner.

For $c,c'\in\bCh(G,N)$, write $c\sim^{\bscl_{G,N}}_a c'$ if $\bscl_{G,N}(c-c')\le a$, where $a\ge0$. Let 
\begin{equation} \label{eq:emb}
\emb \colon [G,N] \to \bCh(G,N), \qquad x \mapsto x
\end{equation}
be the natural map. The following lemma might indicate the reason why the ``coarse geometry of $([G,N],d_{\mathrm{scl}_{G,N}})$'' behaves well.

\begin{lemma}[{compare with \cite[Remark~4.9]{fa}}] \label{lemma coarse homomorphism}
Let $k \geq 1$ and $y_1, \ldots, y_k \in N$. Then the chain
\[
c = y_1 \cdots y_k - (y_1 + \cdots + y_k)
\]
belongs to $\CZ(G,N)$, and $\scl_{G,N}(c) \leq 2^{-1}(k-1)$. In particular, for $x,y \in [G,N]$, we have
\[\emb (xy) \sim^{\bscl_{G,N}}_{2^{-1}} \emb (x) +\emb (y).
\]
\end{lemma}
\begin{proof}
We have
\[
c = -\sum_{i=1}^{k-1} \partial\bigl((y_1 \cdots y_i,\, y_{i+1})\bigr),
\]
and each $(y_1 \cdots y_i, y_{i+1})$ lies in $N \times N$. Hence $c \in \CCC_1(N;\ZZ) \cap \BBB'_1(G,N;\ZZ) = \CZ(G,N)$. For every $\mu \in \QQQ(N)^G$, we have
\[
|\mu(c)| \leq \sum_{i=1}^{k-1} \bigl|\mu(y_1 \cdots y_{i+1}) - \mu(y_1 \cdots y_i) - \mu(y_{i+1})\bigr| \leq (k-1) D(\mu).
\]
Thus Theorem~\ref{theorem gBavard} implies that $\scl_{G,N}(c) \leq 2^{-1}(k-1)$. The last assertion is the case $k = 2$.
\end{proof}

\subsection{Coarse spaces}
For an introduction to coarse spaces, we refer the reader to \cite{Roe}. For a set $X$, let $2^X$ denote the power set of $X$. For $E, F \in 2^{X \times X}$, set
\[ E \circ F = \{ (x,y) \in X \times X \mid \text{there exists $z \in X$ such that $(x,z) \in E$ and $(z,y) \in F$}\},\]
\[ E^{-1} = \{ (x,y) \in X \times X \mid (y,x) \in E\}.\]

\begin{definition}
Let $X$ be a set. A \emph{coarse structure on $X$} is a subset $\EE$ of $2^{X \times X}$ that satisfies the following properties:
\begin{enumerate}[$(1)$]
\item The diagonal set $\Delta_X$  belongs to  $\EE$.

\item If $E \in \EE$ and $F \subset E$, then $F \in \EE$.

\item If $E, F \in \EE$, then $E \circ F \in \EE$.

\item If $E, F \in \EE$, then $E \cup F \in \EE$.

\item If $E \in \EE$, then $E^{-1} \in \EE$.
\end{enumerate}
A \emph{coarse space} is a pair $(X, \EE_X)$ consisting of a set $X$ with a coarse structure $\EE_X$ on $X$.
\end{definition}

\begin{example} \label{example bdd coarse str}
Let $(X,d)$ be a metric space. For $R \ge 0$, set
\[ E_R = \{ (x,y) \in X \times X \; | \; d(x,y) \le R \} = d^{-1} ([0, R]),\]
and set
\[ \EE_d = \{ E \in 2^{X \times X} \; | \; \text{there exists $R \ge 0$ such that $E \subset E_R$}\}.\]
Then $\EE_d$ is a coarse structure on $X$, called the \emph{bounded coarse structure}.
\end{example}

\begin{example}  \label{example almost metric}
  We call a function $d \colon X\times X\to\RR_{\ge 0}$  an \emph{almost metric} if $d(x,x)=0$, $d(x,y)=d(y,x)$, and
there exists a constant $C>0$ such that
\[
d(x,z)\le d(x,y)+d(y,z)+C
\]
for all $x,y,z\in X$. 
 An almost metric  space $(X,d)$ also defines the bounded coarse structure $\EE_d$ on $X$ in the same way as Example~\ref{example bdd coarse str}. 
 
 Fix a constant $C$ as in the definition of an almost metric, and define $d^+ \colon X\times X\to\RR_{\ge 0}$ by 
 \[
d^+(x,y)=
\begin{cases}
d(x,y)+C & \text{if } x\ne y,\\
0 & \text{if } x=y.
\end{cases}
\]
This modification $d^+$ is a genuine metric on $X$ and preserves the bounded coarse structure: $\EE_d = \EE_{d^+}$. In particular, $d^+$ depends on the choice of $C$, but $\EE_{d^+}$ does not.
\end{example}

Let $(X, \EE)$ be a coarse space. For $x,y \in X$ and $E \in \EE$, we write $x \bsim_E\, y$ to mean $(x,y) \in E$. Then, it is straightforward to see that the following hold:
\begin{enumerate}[$(1)$]
\item If $x \bsim_E\, y$, then $y \bsim_{E^{-1}}\, x$.

\item If $x \bsim_E\, y$ and $y \bsim_F\, z$, then $x \bsim_{E \circ F}\, z$.
\end{enumerate}

\begin{definition}
Let $(X, \EE_X)$ and $(Y, \EE_Y)$ be coarse spaces. A map $f \colon X \to Y$ is said to be \emph{bornologous} if $E \in \EE_X$ implies $(f \times f)(E) \in \EE_Y$.
\end{definition}

\begin{definition}
Let $(X, \EE_X)$ be a coarse space, $Y$ a set, and $f$ and $g$ maps from $Y$ to $X$. If
\[ (f \times g)(\Delta_Y) = \{ (f(y), g(y)) \; | \; y \in Y \} \]
 belongs to $\EE_X$, then  we say that $f$ and $g$ are \emph{close}, and write $f \bbsim g$.
\end{definition}

Closeness is an equivalence relation on the set of maps from $Y$ to $X$. Every map close to a bornologous map is bornologous.

\begin{definition}
Let $(X, \EE_X)$ and $(Y, \EE_Y)$ be coarse spaces. 
A \emph{coarse map} is an equivalence class of bornologous maps under the closeness relation $\bbsim$.
\end{definition}

For $E \in 2^{X \times X}$ and for $A \in 2^X$, define
\[ E[A] = E_A = \{ x \in X \; | \; \text{there exists $a \in A$ such that $(x,a) \in E$}\}.\]
For $x \in X$, we write $E_x$ instead of $E_{\{ x\}}$.

\begin{definition}
Let $f \colon (X, \EE_X) \to (Y, \EE_Y)$ be a bornologous map.
\begin{enumerate}[$(1)$]
\item $f$ is called a \emph{coarse embedding} if $F \in \EE_Y$ implies $(f \times f)^{-1}(F) \in \EE_X$.

\item $f$ is said to be \emph{coarsely surjective} if there exists $F \in \EE_Y$ such that $F_{f(X)} = Y$.

\item $f$ is called a \emph{coarse equivalence} if there exists a bornologous map $g \colon (Y, \EE_Y) \to (X, \EE_X)$ such that $g \circ f \bbsim {\rm id}_X$ and $f \circ g \bbsim {\rm id}_Y$.
\end{enumerate}
\end{definition}

The properties of being a coarse embedding, coarsely surjective and a coarse equivalence are invariant up to the closeness relation. The following theorem is well-known:

\begin{theorem}[{\cite{Roe}}]
Let $f \colon (X, \EE_X) \to (Y, \EE_Y)$ be a bornologous map. Then $f$ is a coarse equivalence if and only if $f$ is a coarse embedding and is coarsely surjective.
\end{theorem}

\begin{definition}
Let $(X, \EE_X)$ and $(Y, \EE_Y)$ be coarse spaces. We define the coarse structure $\EE_X \otimes \EE_Y$ on the product set $X \times Y$ as follows:
\[ 
\EE_X \otimes \EE_Y = \{ E \subset (X \times Y) \times (X \times Y) \mid (\pi_1, \pi_3)(E) \in \EE_X, (\pi_2, \pi_4)(E) \in \EE_Y \}, 
\]
where $\pi_i$ denotes the $i$-th projection. 
\end{definition}

For every $E \in 2^{X \times X}$ and $F \in 2^{Y \times Y}$, let $E \otimes F$ be defined by
\[ 
E \otimes F = \{ (x_0, y_0, x_1, y_1) \in X \times Y \times X \times Y \mid (x_0, x_1) \in E, (y_0, y_1) \in F\}. 
\]
It is straightforward to see that $\EE_X \otimes \EE_Y$ is the smallest coarse structure containing the set $\{ E \otimes F \mid E \in \EE_X, F \in \EE_Y \}$.

\begin{definition}
Let $(X, \EE)$ be a coarse space and $A, B \in 2^X$. We say that $A$ is \emph{coarsely included} in $B$, denoted by $A \preccurlyeq B$, if there exists an $E \in \EE$ such that $A \subset E_B$.  If both $A \preccurlyeq B$ and $B \preccurlyeq A$ hold, we say that $A$ and $B$ are \emph{asymptotic}, and we write $A \asymp B$. 
\end{definition}

The relation $\preccurlyeq$ is reflexive and transitive, and the relation $\asymp$ defines an equivalence relation.

\begin{definition}[coarse preimage]
Let $f \colon (X, \EE_X) \to (Y, \EE_Y)$ be a bornologous map, and let $B \in 2^Y$ and $A \in 2^X$. We say that $A$ is a coarse preimage of $B$ with respect to $f$ if both of the following properties are satisfied:
\begin{enumerate}[$(1)$]
\item $f(A) \preccurlyeq B$.

\item For $A' \in 2^X$, $f(A') \preccurlyeq B$ implies $A' \preccurlyeq A$.
\end{enumerate}
\end{definition}

While a coarse preimage may not always exist, it is unique up to asymptoticity whenever it does.

\subsection{Coarse groups}

In this subsection, we review the basic concepts of coarse groups. In our Part~I paper \cite{fa}, we treat groups equipped with bi-invariant metrics as main examples of coarse groups. In the present paper, we go into the general setting of the theory of coarse groups. In Leitner--Vigolo \cite{LV}, a coarse group is defined as a group object in the category of coarse spaces and coarse maps, i.e., the  multiplication and inversion are coarse maps, and the group
axioms hold in the coarse category.  In this paper, we formulate a coarse group as a coarse space equipped with a binary operation that is bornologous.  This is a slight modification intended to simplify the subsequent discussion in this paper and does not constitute a substantial change.

\begin{definition}
A \emph{coarse group} $(G, m, e_G, c)$ consists of the following data: a coarse space $G = (G, \EE_G)$, a bornologous map $m \colon (G \times G, \EE_G \otimes \EE_G) \to (G, \EE_G)$, an element $e_G$ of $G$ called the \emph{unit of $G$}, and a bornologous map $c \colon G \to G$. We require that the following diagrams are commutative up to closeness:
\[ \xymatrix{
G \times G \times G \ar[r]^-{m \times {\rm id}_G} \ar[d]_-{{\rm id}_G \times m} & G \times G \ar[d]^-{m} \\
G \times G \ar[r]^-{m} & G,
} \quad \xymatrix{
G \ar[r]^-{(e_G, {\rm id}_G)} \ar[rd]_-{{\rm id}_G} & G \times G \ar[d]^-{m} & \ar[l]_-{({\rm id}_G, e_G)} \ar[ld]^-{{\rm id}_G} G \\
& G &
} \quad \xymatrix{
G \ar[r]^-{({\rm id}_G, c)} \ar[d]_-{(c, {\rm id}_G)} \ar[rd]^{e_G} & G \times G \ar[d]^{m} \\
G \times G \ar[r]^m & G 
}\]
Here, $e_G \colon G \to G$ means the constant map at $e_G$. When $G$ is coarsely abelian, we write $0_G$ for $e_G$, and we use the same convention for the constant map. We often say ``$G$ is a coarse group'' instead of ``$(G, m, e_G, c)$ is a coarse group.''
\end{definition}

We abbreviate $m(x,y)$ to $xy$  and write $x^{-1}$ for $c(x)$. 
Henceforth, for elements $x, y, z$ of a coarse group, the expression $xyz$ is understood as $((xy)z)$. Similarly, $xyzw$ and longer products are regarded in the same manner.

\begin{example}
Let $d$ be a bi-invariant metric on a group $G$. Then $G$ forms a coarse group under the coarse structure induced by $d$. In particular, we regard $([G,N], d^+_{\scl_{G,N}})$ as a coarse group in this way.
\end{example}

\begin{remark} \label{remark d_scl d_scl+}
Since $\scl_{G,N}(xy) \leq \scl_{G,N}(x) + \scl_{G,N}(y)+\frac{1}{2}$, $\dGN$ is an almost metric on $[G,N]$. 
By modifying $\dGN$ as in Example~\ref{example almost metric} with $C=1/2$, we obtain a bi-invariant metric $\dGN^+$ on $[G,N]$; 
however, this modification is introduced merely to ensure the metric axioms, and since $\mathcal{E}_{\dGN}=\mathcal{E}_{\dGN^+}$, it does not affect the coarse structure. Therefore, $([G,N], \dGN)$ and $([G,N], \dGN^+)$ define the same coarse group, and we shall work with the more natural $([G,N],\dGN)$ throughout the paper.
\end{remark}

\begin{definition} \label{coarse homomorphism}
Let $G$ and $H$ be coarse groups. A bornologous map $f \colon G \to H$ is called a \emph{coarse homomorphism} if there exists $E \in \EE_H$ such that $(f(xy), f(x)f(y)) \in E$ for every $x,y \in G$. A \emph{coarse isomorphism} is a coarse homomorphism between coarse groups that is coarsely surjective and a coarse embedding.
\end{definition}

\begin{example} \label{example natural embedding}
Let $\emb \colon [G,N] \to \bCh(G,N)$ be the natural map  in \eqref{eq:emb}.  Then $\emb$ is a coarse homomorphism from $([G,N], \dGN)$ to $(\bCh(G,N), \bscl_{G,N})$ by Lemma~\ref{lemma coarse homomorphism}. Note that $\emb$ is also a coarse embedding; $\dGN(x,y) \sim_{2^{-1}} \bscl_{G,N}(\emb(x) - \emb(y))$ for $x,y \in [G,N]$.
\end{example}

The following lemma is immediately deduced from the definition of coarse groups, and we will use it frequently in this paper without explicit reference.

\begin{lemma}[{\cite{LV}}] \label{lemma uniform}
Let $G$ be a coarse group and let $E \in \EE_G$. Then the following hold:
\begin{enumerate}[$(1)$]
\item There exists $F \in \EE_G$ such that $(xz, yz) \in F$ for all $(x,y) \in E$ and $z \in G$.

\item There exists $F \in \EE_G$ such that $(zx, zy) \in F$ for all $(x,y) \in E$ and $z \in G$.
\end{enumerate}
\end{lemma}

\begin{definition}
Let $G$ be a coarse group. A subset $H$ of $G$ is called a \emph{coarse subgroup} if it satisfies the following properties:
\begin{itemize}
\item $\{ e_G\} \preccurlyeq H$

\item $m(H \times H) \preccurlyeq H$

\item $H^{-1} \preccurlyeq H$
\end{itemize}
\end{definition}

By \cite[Proposition~5.3.2]{LV}, a coarse subgroup $H\subset G$, equipped with the subspace coarse structure, admits a coarse group structure for which the inclusion $H\hookrightarrow G$ is a coarse homomorphism. This structure is unique up to closeness of the operations. We shall always regard coarse subgroups as coarse groups in this way.

\begin{definition}
Let $G$ be a coarse group and let $H$ be a coarse subgroup of $G$. $H$ is called a \emph{coarse normal subgroup} if 
 $\{ghg^{-1}\mid g\in G,\ h\in H\}\preccurlyeq H$. 
\end{definition}

\begin{definition}[coarse kernel] \label{definition coarse kernel}
Let $f \colon G \to H$ be a coarse homomorphism between coarse groups. If a coarse preimage of $\{e_H\}$ with respect to $f$ exists, it is called the \emph{coarse kernel} of $f$ and is denoted by $\KKer(f)$.
\end{definition}

A coarse kernel $K \subset G$ is a coarse subgroup of $G$ (see \cite[Lemma~7.1.1]{LV} for example).
Here, we equip $K$ with the coarse group structure described above.

\begin{remark} \label{remark coarse kernel}
Let $j \colon K \to G$ and $f \colon G \to H$ be coarse homomorphisms. Then it is straightforward to see that the following conditions are equivalent:
\begin{enumerate}[(1)]
\item $j$ is a coarse embedding and $j(K)$ is the coarse kernel of $f \colon G \to H$ in the sense of Definition~\ref{definition coarse kernel}.

\item $f \circ j \bbsim e_H$, and for every coarse space $X$ and every bornologous map $g \colon X \to G$ such that $f \circ g \bbsim e_H$, there exists a bornologous map $g' \colon X \to K$, which is unique up to closeness, such that $j g' \bbsim g$.
\end{enumerate}
In this paper, we say that $j \colon K \to G$ is the coarse kernel of $f$ if $j$ satisfies one of these conditions.
\end{remark}

\begin{definition} \label{definition exact}
Let $G_0$, $G_1$ and $G_2$ be coarse groups, and let $i \colon G_0 \to G_1$ and $p \colon G_1 \to G_2$ be coarse homomorphisms. We say that the sequence
\[ 1 \to G_0 \xrightarrow{i} G_1 \xrightarrow{p} G_2 \]
is \emph{exact} if $i \colon G_0 \to G_1$ is the coarse kernel of $p$ in the sense of Remark~\textup{\ref{remark coarse kernel}}. If moreover $p$ is coarsely surjective, we call
\[ 1 \to G_0 \xrightarrow{i} G_1 \xrightarrow{p} G_2 \to 1 \]
a \emph{short exact sequence}. For coarsely abelian groups, we write $0$ instead of $1$.
\end{definition}

\begin{definition}[{\cite[Definition~6.3.3]{LV}}]
A coarse group $(G, \EE_G)$ is a \emph{coarsely abelian group} if the following equivalent conditions hold:
\begin{enumerate}[$(1)$]
\item Define $\tau \colon G \times G \to G \times G$ by $\tau (x,y) = (y,x)$. Then the diagram
\[ \xymatrix{
G \times G \ar[r]^-{\tau} \ar[rd]^-{m} & G \times G \ar[d]^-m \\
 & G
}\]
commutes up to closeness.

\item $\{ (xy, yx) \mid x,y \in G \} \in \EE_G$.

\item $\{ (gxg^{-1}, x) \mid g,x\in G\} \in \EE_G$.

\item $\{ xyx^{-1}y^{-1}  \mid x,y\in G  \} \asymp \{ e_G \}$.

\item $m \colon G \times G \to G$ is a coarse homomorphism.

\item $c \colon G \to G$ given by $x \mapsto x^{-1}$ is a coarse homomorphism.
\end{enumerate}
\end{definition}

The coarse group $([G,N], \dGN)$ is coarsely abelian. This fact can be easily checked, or follows from Example~\ref{example natural embedding}.

It is straightforward to check the following lemma from the definition.

\begin{lemma} \label{lemma homomorphism sum}
Let $G$ be a coarse group and $H$ a coarsely abelian group. Suppose that $f, g \colon G \to H$ are coarse homomorphisms. Then the maps $-f$ and $f+g$ defined by
\[ (-f)(x) = -f(x), \quad (f + g)(x) = f(x) + g(x)\]
are also coarse homomorphisms.
\end{lemma}

We will use the following theorem in the proof of Proposition~\ref{proposition splitting}.

\begin{theorem}[{\cite[Proposition~1.4.6]{LV}}] \label{theorem linear}
Let $f \colon (\RR^m, \| \cdot \|_1) \to (\RR^n, \| \cdot \|_1)$ be a coarse homomorphism. Then there exists a unique linear map $T \colon \RR^m \to \RR^n$ such that $f$ and $T$ are close.
\end{theorem}

\section{Duality} \label{section functional analysis}

Let $N$ and $L$ be normal subgroups of $G$ such that
$N\leqslant L$. The natural map
\[
T_{N,L}\colon \Ch(G,N)\longrightarrow\Ch(G,L),
\qquad [c]\longmapsto[c],
\]
is norm-nonincreasing. It therefore extends uniquely to a bounded
linear map
\[
\bT_{N,L}\colon\bCh(G,N)\longrightarrow\bCh(G,L).
\]

As explained in Section~\ref{section introduction}, the kernel $\Ker(\bT_{N,L})$ plays an important role in the proof of Theorem~\ref{theorem A}. 
The purpose of this section is to identify its continuous dual in terms of $\WW(G,L,N)$ under a certain condition. In this section, $N \leqslant L$ are arbitrary normal subgroups of an arbitrary group $G$. Recall from Subsection~\ref{subsec=qm} that $\WW(G,L,N)$ is the cokernel of $\widehat{i^*}_{N,L}\colon \QQQ(L)^G/\HHH^1(L)^G\to\QQQ(N)^G/\HHH^1(N)^G$. Our only assumption is that the image of $\widehat{i^*}_{N,L}$ is norm-closed; this holds whenever $\WW(G,L,N)$ is finite-dimensional (Lemma~\ref{lemma finite-codimensional}). Under this assumption, $\WW(G,L,N)$ carries the quotient norm $\bD_N$; we use $2\bD_N$, the factor $2$ matching that in Bavard duality. We use the notion of a duality pairing introduced in Definition~\ref{def=dualitypairing}.

\begin{theorem} \label{theorem W duality}
  Assume that the image of $\widehat{i^*}_{N,L}$ is norm-closed.
  Then, the following map
\[ 
(\WW(G,L,N), 2\bD_N) \times (\Ker(\bT_{N,L}), \bscl_{G,N})\to \RR;\quad ([\mu],\alpha)\mapsto \mu(\alpha)
\]
is a duality pairing for $(\Ker(\bT_{N,L}), \bscl_{G,N})$, and thus  $(\Ker(\bT_{N,L}), \bscl_{G,N})^*$ is isometrically isomorphic to $(\WW(G,L,N), 2\bD_N)$.  In particular, if $\WW(G,L,N)$ is finite-dimensional, then the image of $\widehat{i^*}_{N,L}$ is norm-closed, and hence the conclusion above holds.
\end{theorem}

 The case $L=G$ was proved in \cite{fa}, and the same argument
applies to the present setting. For completeness, we include the
proof. 
For the proof of the last assertion, we will employ the following general fact in functional analysis.

\begin{lemma}[{see \cite[Corollary~3.8]{fa} for example}] \label{lemma finite-codimensional}
Let $X$ and $Y$ be real Banach spaces. Let $S \colon X \to Y$ be a bounded linear operator. If the quotient space $Y / S(X)$ is finite-dimensional, then $S(X)$ is norm-closed in $Y$.
\end{lemma}

 We use the convention for continuous duals and adjoint operators introduced after Definition~\ref{def=dualitypairing}. 

\begin{proof}[Proof of Theorem~\textup{\ref{theorem W duality}}]
By Theorem~\ref{theorem duality}, $(\QQQ(N)^G / \HHH^1(N)^G,2\DDD_N)$ is identified with $(\bCh(G,N),\bscl_{\Gg,\Ng})^*$ via the duality pairing
\[
\langle \cdot,\cdot\rangle_{G,N}\colon (\QQQ(N)^G / \HHH^1(N)^G)\times \bCh(G,N)  \to \RR  ;\quad ([\mu],\alpha)\mapsto\mu(\alpha)
\]
for $(\bCh(G,N),\bscl_{\Gg,\Ng})$.
In the same manner, $(\QQQ(L)^G / \HHH^1(L)^G,2\DDD_L)$ is identified with $(\bCh(G,L),\bscl_{\Gg,L})^*$ via the duality pairing
\[
\langle \cdot,\cdot\rangle_{G,L}\colon (\QQQ(L)^G / \HHH^1(L)^G)\times \bCh(G,L)  \to \RR ;\quad ([\phi],\alpha)\mapsto\phi(\alpha)
\]
for $(\bCh(G,L),\bscl_{\Gg,L})$.
For every $\phi\in \QQQ(L)^G$ and $c\in \CR(G,N)$, we have
\[
\langle  \widehat{i^*}_{N,L} [\phi],[c]\rangle_{G,N}=\phi(c)=\langle [\phi],\bT_{N,L}[c]\rangle_{G,L}.
\]

Since $\Ch(G,N)$ is dense in $\bCh(G,N)$ and both sides are continuous in the second variable, this identity extends to every $\alpha\in\bCh(G,N)$. Hence $\bT_{N,L}^{\dagger}=\widehat{i^*}_{N,L}$ under the identifications above.

Since \(\bT_{N,L}\) is the bounded extension of the norm-nonincreasing map \(T_{N,L}\), 
the operator norm of $\bT_{N,L}$ is at most $1$. 
The spaces $(\bCh(G,N),\bscl_{\Gg,N})$ and $(\bCh(G,L),\bscl_{\Gg,L})$ are Banach spaces by construction. 
By hypothesis, the image of $\widehat{i^*}_{N,L}$ is norm-closed.
Applying Corollary~\ref{cor=Banach} to $S=\bT_{N,L}$ and using $\bT_{N,L}^{\dagger}=\widehat{i^*}_{N,L}$, we obtain the desired duality pairing.

The dualities above show that
$\bigl(\QQQ(L)^G/\HHH^1(L)^G,2D_L\bigr)$
and
$\bigl(\QQQ(N)^G/\HHH^1(N)^G,2D_N\bigr)$
are Banach spaces. If $\WW(G,L,N)$ is finite-dimensional, then
$\widehat{i^*}_{N,L}$ has finite-codimensional range. Moreover,
$\widehat{i^*}_{N,L}$ is bounded because it is the adjoint of
$\bT_{N,L}$. Lemma~\ref{lemma finite-codimensional} therefore
implies that its range is norm-closed. This proves the last
assertion.

\end{proof}

\begin{corollary} \label{corollary finite-dimensional}
 Assume that $\WW(G,L,N)$ is finite-dimensional. Then, we have an isometric isomorphism $(\Ker(\bT_{N,L}),\bscl_{G,N})\cong (\WW(G,L,N),2\bD_N)^*$. In particular, the dimension of $\Ker(\bT_{N,L})$ equals that of $\WW(G,L,N)$. Moreover, $\Imm(\bT_{N,L})$ is norm-closed.
\end{corollary}

\begin{proof}
By Theorem~\ref{theorem W duality}, we have an isometric isomorphism $\Ker(\bT_{N,L}) ^{**}\cong \WW(G,L,N)^*$. Since $\Ker(\bT_{N,L})$ is naturally embedded in $\Ker(\bT_{N,L})^{**}$, this implies that $\Ker(\bT_{N,L})$ is finite-dimensional. Hence $\Ker(\bT_{N,L})$ is reflexive, and  we have the desired isometric isomorphism. Finally, $\bT_{N,L}^{\dagger}=\widehat{i^*}_{N,L}$ (see the proof of Theorem~\ref{theorem W duality}) has norm-closed image by Theorem~\ref{theorem W duality}, and hence $\Imm(\bT_{N,L})$ is norm-closed by Theorem~\ref{theorem CR}.
\end{proof}

\section{Elements of $\Ker(\bT_{N,L})$} \label{section kernel element}

Let $N$ and $L$ be normal subgroups of $G$ such that $N \leqslant L$ and $[G,L] \leqslant N$. For $g \in G$ and $x \in L$, consider the sequence
\[ \left( \frac{[g,x^n]}{n} \right)_{n \in \NN}\]
of chains. This sequence is not necessarily contained in $\CQ(G,N)$. However, as shown in the following lemma, it behaves like a Cauchy sequence with respect to $\scl_{G,N}$. Apart from $[G,L] \leqslant N$, no assumption on $G$, $L$ and $N$ (such as the lower central series or finite generation of $G$) is used in this section. Throughout this section, for $c_1, c_2 \in \bCh(G,N)$ and $a \ge 0$, we write $c_1 \sim^{\bscl_{G,N}}_a c_2$ to mean $\bscl_{G,N}(c_1 - c_2) \le a$.

\begin{lemma} \label{lemma Cauchy}
Let $N$ and $L$ be normal subgroups of $G$ such that $N \leqslant L$ and $[G,L] \leqslant N$. For $g \in G$ and $x \in L$, the following hold:
\begin{enumerate}[$(1)$]
\item For positive integers $m$ and $n$, we have
\[ \frac{[g,x^m]}{m} - \frac{[g,x^n]}{n} \in \CQ(G,N).\]
\item $\displaystyle \lim_{m,n \to \infty} \scl_{G,N} \left( \frac{[g,x^m]}{m} - \frac{[g,x^n]}{n} \right) = 0$.

\end{enumerate}
\end{lemma}

\begin{proof}
For a positive integer $k$, set $a_k=[g,x^k]$. Since
$x\in L$ and $a_1=[g,x]\in N$, the elements $x$ and $a_1$
commute in $G/[G,N]$, and hence $a_k=a_1^k$ in $G/[G,N]$.
In particular, $a_m^n=a_n^m$ in $G/[G,N]$. By
Lemma~\ref{lem=domainG}, this implies
\[
\frac{a_m}{m}-\frac{a_n}{n}\in\CQ(G,N).
\]
This proves \textup{(1)}.

We now prove \textup{(2)}. The commutator identity gives
\[
a_{mn}
=
a_m\cdot x^m a_m x^{-m}\cdots
x^{m(n-1)}a_mx^{-m(n-1)}.
\]
The argument above also shows that
\[
n a_m-a_{mn}\in\CZ(G,N)
\quad\text{and}\quad
a_{mn}-m a_n\in\CZ(G,N).
\]
Thus, the rational chains appearing below belong to
$\CQ(G,N)$.

Since every $\mu\in\QQQ(N)^G$ is $G$-invariant, we have
\[
|\mu(a_{mn})-n\mu(a_m)|
\leq (n-1)D(\mu).
\]
Similarly,
\[
|\mu(a_{mn})-m\mu(a_n)|
\leq (m-1)D(\mu).
\]
The generalized mixed Bavard duality therefore gives
\begin{align} \label{scl cauchy}
\scl_{G,N}\left(
 \frac{a_m}{m}-\frac{a_{mn}}{mn}
\right)
\leq\frac{n-1}{2mn}
\end{align}
and
\[
\scl_{G,N}\left(
 \frac{a_{mn}}{mn}-\frac{a_n}{n}
\right)
\leq\frac{m-1}{2mn}.
\]
Thus, by the triangle inequality,
\[
\scl_{G,N}\left(
 \frac{a_m}{m}-\frac{a_n}{n}
\right)
\leq
\frac{n-1}{2mn}+\frac{m-1}{2mn}
=
\frac{m+n-2}{2mn}.
\]
Since the right-hand side tends to $0$ as $m,n\to\infty$,
this proves \textup{(2)}.
\end{proof}

\begin{proposition} \label{proposition kernel element}
Let $N$ and $L$ be normal subgroups of $G$ such that $N \leqslant L$ and $[G,L] \leqslant N$. Let $g_1, \cdots, g_k \in G$ and $x_1, \cdots, x_k \in L$ such that
\[ [g_1, x_1] \cdots [g_k, x_k] \in [G,N].\]
Then the sequence 
\[ (c_n)_{n \in \NN} = \left( \frac{[g_1, x_1^n]}{n} + \cdots + \frac{[g_k, x_k^n]}{n}\right)_{n \in \NN} \]
is a Cauchy sequence in $(\bCh(G,N), \bscl_{G,N})$, and its norm-limit  belongs to  $\Ker(\bT_{N,L})$.
\end{proposition}

\begin{proof}
We first show that every $c_n$ is contained in $\CQ(G,N)$. By Lemma~\ref{lem=domainG}, it suffices to show that $[g_1, x_1^n] \cdots [g_k, x_k^n] \in [G,N]$. As in the proof of Lemma~\ref{lemma Cauchy}, we have $[g_i, x_i^n] = [g_i, x_i]^n$ in $G/[G,N]$. Since $[g_i, x_i], [g_j, x_j] \in N$ commute in $G/[G,N]$, we conclude that in $G/[G,N]$,
\[ [g_1, x_1^n] \cdots [g_k, x_k^n] = [g_1, x_1]^n \cdots [g_k, x_k]^n = ([g_1, x_1] \cdots [g_k, x_k])^n.\]
Since $[g_1, x_1] \cdots [g_k, x_k] \in [G,N]$, we have $[g_1, x_1^n] \cdots [g_k, x_k^n] \in [G,N]$. Thus we have proven that $(c_n)_{n \in \NN}$ is a sequence in $\CQ(G,N)$, and hence is in $\bCh(G,N)$. Then, it immediately follows from Lemma~\ref{lemma Cauchy} (2) that $(c_n)_{n \in \NN}$ is a Cauchy sequence in $\bCh(G,N)$.

Let $\alpha$ be the norm-limit of $(c_n)_{n \in \NN}$ in $\bCh(G,N)$. 
 Under the map $\bT_{N,L}$, the element $\alpha$ is represented by
the same sequence $(c_n)$. On the other hand,
\[
\bscl_{G,L}(c_n)\leq\frac{k}{2n},
\]
because each $[g_i,x_i^n]$ is a single $(G,L)$-commutator.
Therefore $(c_n)$ converges to zero in $\bCh(G,L)$, and hence
$\bT_{N,L}(\alpha)=0$. 
\end{proof}

Here, while the norm-limit $\alpha$ of $(c_n)_{n \in \NN}$ in Proposition~\ref{proposition kernel element} is $0$ in $\bCh(G,L)$, it is not necessarily $0$ in $\bCh(G,N)$. In the next section, we demonstrate the existence of such an example for the case where $N =  \gamma_{q+1} (G)$ and $L =  \gamma_{q} (G)$. Here $\gamma_q(G)$ denotes the $q$-th term of the lower central series (see Subsection~\ref{subsection notation}).

The next goal is Proposition~\ref{proposition group lifting}, which shows that  certain linear maps into $\Ker(\bT_{N,L})$ can be realized, up to closeness in the chain space, by coarse homomorphisms from $(\ZZ^l,\|\cdot\|_1)$ to $[G,N]$. 

\begin{lemma} \label{lemma line lifting}
Let $N, L, g_i, x_i$, and $(c_n)_{n \in \NN}$ be as in Proposition~\textup{\ref{proposition kernel element}}, and let $\alpha$ be the norm-limit of $(c_n)_{n \in \NN}$ in $\bCh(G,N)$. Then, for every $m \in \ZZ$, we have
\[ [g_1, x_1^m] + \cdots + [g_k, x_k^m] \sim^{\bscl_{G,N}}_{2^{-1}k} m \alpha.\]
\end{lemma}

\begin{proof}
Since the case where $m = 0$ is trivial, we  first assume $m > 0$.  By Lemma~\ref{lemma Cauchy}, we have $\dfrac{[g_i, x_i^{mn}]}{n} - [g_i, x_i^m] \in \CR(G,N)$. 
 Multiplying
\eqref{scl cauchy} by $m$, with $g=g_i$ and $x=x_i$, gives
\[
\scl_{G,N}\left(
 \frac{[g_i,x_i^{mn}]}{n}-[g_i,x_i^m]
\right)
\leq\frac{n-1}{2n}
\leq\frac12.
\] 
Therefore, we obtain
\[ \scl_{G,N} \left( m \cdot \Big(\frac{[g_1, x_1^{mn}]}{mn} + \cdots + \frac{[g_k, x_k^{mn}]}{mn}\Big) - \big([g_1, x_1^m] + \cdots + [g_k, x_k^m] \big)\right) \le \frac{k}{2}. \]
Taking the limit as $n\to\infty$, we obtain
\[ \scl_{G,N} \Big(m\alpha - \big( [g_1, x_1^m] + \cdots + [g_k, x_k^m] \big) \Big) \le \frac{k}{2}. \]
 For $m<0$, the data $g_i,x_i^{-1}$ satisfy the same assumptions as $g_i,x_i$. Hence we may apply the case $m>0$ to $x_i^{-1}$ and $-m$. Since $[g_i,x_i^{-n}]$ is conjugate to $[g_i,x_i^n]^{-1}$, the corresponding norm-limit is $-\alpha$. Hence the same estimate holds for $m<0$. 
\end{proof}

\begin{lemma} \label{lemma lifting}
Let $N$ and $L$ be normal subgroups of $G$ such that $N \leqslant L$ and $[G,L] \leqslant N$. Let $g_1, \cdots, g_k \in G$ and $x_1, \cdots, x_k \in L$ such that
\[ [g_1, x_1] \cdots [g_k, x_k] \in [G,N].\]
Set
\[ \alpha = \nlim{\bscl_{G,N}}_{n \to \infty} \left( \frac{[g_1, x_1^n]}{n} + \cdots + \frac{[g_{k}, x_{k}^n]}{n}\right) \in \bCh(G,N). \]
For $m \in \ZZ$, set $\alpha^{(m)} = [g_1, x_1^m] \cdots [g_k, x_k^m] \in [G,N]$. Then $\alpha^{(m)} \sim^{\bscl_{G,N}}_{2^{-1}(2k-1)} m \alpha$.

\end{lemma}
\begin{proof}
Each $[g_i, x_i^m]$ belongs to $N$ since $[G,L] \leqslant N$. By Lemma~\ref{lemma line lifting} and Lemma~\ref{lemma coarse homomorphism} applied to $y_i = [g_i, x_i^m]$, we have
\[ m \alpha \sim^{\bscl_{G,N}}_{2^{-1} k} [g_1, x_1^m] + \cdots + [g_k, x_k^m] \sim^{\bscl_{G,N}}_{2^{-1}(k-1)} [g_1, x_1^m]  \cdots  [g_k, x_k^m] = \alpha^{(m)}.\qedhere\]
\end{proof}

\begin{proposition} \label{proposition group lifting}
Let $N$ and $L$ be normal subgroups of $G$ such that $N \leqslant L$ and $[G,L] \leqslant N$
 and let $l$ be a non-negative
integer. For each $i\in\{1,\dots,l\}$, choose data satisfying the hypotheses of Lemma~\textup{\ref{lemma lifting}}, with $k=k_i$, and denote the resulting elements by
\[
\alpha_i\in\bCh(G,N)
\quad\text{and}\quad
\alpha_i^{(m)}\in[G,N]\qquad(m\in\ZZ).
\] 
Define $\Phi \colon (\ZZ^l, \| \cdot \|_1) \to [G,N]$ by
\[ \Phi(m_1, \dots, m_l) = \alpha_1^{(m_1)} \cdots \alpha_l^{(m_l)}. \]
Then, $\Phi$ is a coarse homomorphism, and it satisfies
\[ \emb \circ \Phi(m_1, \dots, m_l) \sim^{\bscl_{G,N}}_{kl} m_1 \alpha_1 + \dots + m_l \alpha_l, \]
where $k = \max \{ k_1, \dots, k_l\}$ (with $k=0$ when $l=0$) and $\emb \colon [G,N] \to \bCh(G,N)$ is the natural map  in \eqref{eq:emb}. 
\end{proposition}

\begin{proof}
 By repeated application of Lemma~\ref{lemma coarse homomorphism}, we have 
\begin{align} \label{a1}
\alpha_1^{(m_1)} \cdots \alpha_l^{(m_l)} \sim^{\bscl_{G,N}}_{2^{-1} l} \alpha_1^{(m_1)} + \cdots + \alpha_l^{(m_l)}.
\end{align}
By Lemma~\ref{lemma lifting}, we have
\begin{align} \label{a2}
\alpha_i^{(m_i)} \sim^{\bscl_{G,N}}_{2^{-1}(2k-1)} m_i \alpha_i.
\end{align}
Combining \eqref{a1} and \eqref{a2}, we obtain
\begin{align} \label{a3}
\emb\circ \Phi(m_1, \dots, m_l) \sim^{\bscl_{G,N}}_{2^{-1}(l+l(2k - 1))} m_1 \alpha_1 + \dots + m_l \alpha_l.
\end{align}

We show that $\Phi$ is a coarse homomorphism. Define $\Phi' \colon \ZZ^l \to \bCh(G,N)$ by
\[ \Phi'(m_1, \dots, m_l) = m_1 \alpha_1 + \dots + m_l \alpha_l. \]

The map $\Phi'$ is a Lipschitz homomorphism and hence is
bornologous. Since \eqref{a3} shows that
$\emb\circ\Phi$ is close to $\Phi'$, the map
$\emb\circ\Phi$ is also bornologous. Since $\emb$ is a coarse
embedding, it follows that $\Phi$ is bornologous.

Moreover, for every $u,v\in\ZZ^l$, we have
\[
\begin{aligned}
\emb\bigl(\Phi(u+v)\bigr)
&\bbsim \Phi'(u+v)\\
&=\Phi'(u)+\Phi'(v)\\
&\bbsim
\emb\bigl(\Phi(u)\bigr)+\emb\bigl(\Phi(v)\bigr)\\
&\bbsim
\emb\bigl(\Phi(u)\Phi(v)\bigr).
\end{aligned}
\]
The last closeness follows from Lemma~\ref{lemma coarse homomorphism}, applied to $x=\Phi(u)$ and $y=\Phi(v)$.
Since $\emb$ is a coarse embedding, we conclude that
\[
\Phi(u+v)\bbsim\Phi(u)\Phi(v).
\]
Thus $\Phi$ is a coarse homomorphism. 
\end{proof}

From the discussion in this section, we are interested in  nonzero elements of $\Ker(\bT_{N,L})$ that can be expressed as the norm-limit  of the form 
\[ \nlim{\bscl_{G,N}}_{n \to \infty} \left( \frac{[g_1, x_1^n]}{n} + \cdots + \frac{[g_k, x_k^n]}{n} \right) \quad (g_i \in G, x_i \in L). \]
In the next section, we will see that  $\Ker(\bT_{N,L})$ admits a basis consisting of such elements when
$\dim_{\RR}\WW(G,L,N)<\infty$, $N =  \gamma_{q+1} (G)$, and $L =  \gamma_{q} (G)$ (where $\gamma_q(G)$ denotes the $q$-th term of the lower central series of $G$).

\section{Core} \label{section core}

In this section, we analyze the coarse kernel of $\bT_{N,L} \colon \bCh(G,N) \to \bCh(G,L)$ 
for the case $ N =  \gamma_{q+1} (G)$ and $L =  \gamma_{q} (G)$, and prove Theorem~\ref{theorem B} and Theorem~\ref{theorem C}. The first subsection only uses $[G,L] \leqslant N$. From Subsection~\ref{subsection core} on, we work in Setting~\ref{setting}, and in Subsection~\ref{subsection 5.3} we additionally assume that $\WW(G,L,N)$ is finite-dimensional. Finite generation of $G$ is used only in the proof of Theorem~\ref{theorem C}, to guarantee this finite-dimensionality (Corollary~\ref{corollary fg nilpotent W}~\textup{(2)}) and the coarse surjectivity (Corollary~\ref{corollary xi equivalence}).

\subsection{Characterization of extendability}

Let $N$ and $L$ be normal subgroups of $G$ such that $N \leqslant L$.
Let $\mu$ be a homogeneous quasimorphism on $N$, not necessarily $G$-invariant.
Recall the defect formula
\begin{equation}\label{eq:defect-commutators}
D(\mu)=\sup_{a,b\in N}|\mu([a,b])|
\end{equation}
(see \cite[Lemma~2.24]{Calegari}  for example).
We consider two related quantities obtained by allowing commutators involving elements outside $N$:
\begin{align*}
\DDD_{G,N}(\mu) & = \sup \{ |\mu([g,x])| \mid g \in G,\ x \in N\}, \\
\DDD_{G,N,L}(\mu) & = \sup \{ |\mu([g,h])| \mid \text{$g \in G$,\ $h \in L$ such that $[g,h] \in N$} \}.
\end{align*}
These quantities may be infinite, and
\[
D(\mu)\leq \DDD_{G,N}(\mu)\leq \DDD_{G,N,L}(\mu).
\]

\begin{lemma} \label{lemma ch inv}
Let $G$ be a group, and $N$ a normal subgroup of $G$. Let $\mu$ be a homogeneous quasimorphism on $N$. Then $\mu$ is $G$-invariant if and only if $\DDD_{G,N}(\mu) < \infty$.
\end{lemma}
\begin{proof}
Suppose that $\mu$ is $G$-invariant. Then we have
\[ \mu([g,x]) \sim_{D(\mu)} \mu(gxg^{-1}) + \mu(x^{-1}) = 0.\]
Hence we have $|\mu([g,x])| \le \DDD(\mu)$. This means $\DDD_{G,N}(\mu) \le \DDD(\mu) < \infty$.

Next suppose $\DDD_{G,N}(\mu) < \infty$. Suppose that there exist $g \in G$ and $x \in N$ such that $\mu(gxg^{-1}) \ne \mu(x)$. Since $gx^ng^{-1}=(gxg^{-1})^n$, we may choose a positive integer $n$ such that, by homogeneity,
\[ |\mu(gx^n g^{-1}) - \mu(x^n)| = n |\mu(gxg^{-1}) - \mu(x)| > \DDD_{G,N}(\mu) + \DDD(\mu).\]
Then we have
\[ \DDD_{G,N}(\mu) + \DDD(\mu) < |\mu(gx^n g^{-1}) - \mu(x^n)| \sim_{\DDD(\mu)} |\mu([g,x^n])| \le \DDD_{G,N}(\mu).\]
This is a contradiction.
\end{proof}

The first part of the proof, together with \eqref{eq:defect-commutators}, shows that $\DDD_{G,N}(\mu)=D(\mu)$ whenever $\mu\in\QQQ(N)^G$.

\begin{proposition} \label{proposition ch ext}
Let $G$ be a group, $L$ and $N$ normal subgroups of $G$ such that $N \leqslant L$, and let $\mu \in \QQQ(N)^G$. Assume that $[G,L] \leqslant N$. Then the following are equivalent:
\begin{enumerate}[$(1)$]
\item There exists $\psi \in \QQQ(L)^G$ such that $\psi|_N = \mu$.

\item $\DDD_{G,N,L}(\mu) < \infty$.
\end{enumerate}
\end{proposition}
\begin{proof}

Suppose that (1) holds, and choose $\psi\in\QQQ(L)^G$ such that $\psi|_N=\mu$.
By the first part of the proof of Lemma~\ref{lemma ch inv}, applied to $\psi$ on $L$, we have
\[
\DDD_{G,N,L}(\mu)=\DDD_{G,L}(\psi)\leq D(\psi)<\infty.
\]
Here the equality follows from $[G,L]\leqslant N$ and $\psi|_N=\mu$.
Thus (2) holds.

Conversely, suppose that (2) holds. Set
\[ \FF = \{ (\psi, H) \mid \text{$N \leqslant H \leqslant L,\ H \trianglelefteqslant G, \psi \in \QQQ(H)^G$ such that $\psi |_N = \mu$} \}.\]
Define a partial order $\preccurlyeq$ on $\FF$ by
\[ (\psi_0, H_0) \preccurlyeq (\psi_1, H_1) \quad \text{if} \quad \text{$H_0 \leqslant H_1$ and $\psi_1|_{H_0} = \psi_0$.}\]
 For every $(\psi,H)\in\FF$, we have, by the defect formula \eqref{eq:defect-commutators} applied to $\psi$ on $H$,
\[
D(\psi)=\sup_{a,b\in H}|\psi([a,b])|
\leq \DDD_{G,N,L}(\mu),
\]
since $[H,H]\leqslant N$ and $\psi|_N=\mu$. Hence 
 every totally ordered subset of $\FF$ has an upper bound. Note that $(\mu, N)$ belongs to $\FF$, and hence that $\FF \ne \emptyset$. Hence, it follows from Zorn's lemma that $(\FF, \preccurlyeq)$ has a maximal element $(\psi, H)$.

It suffices to show $H = L$. Suppose that $H \subsetneq L$, and let $x \in L \setminus H$. Let $K$ be the smallest normal subgroup of $G$ containing $H$ and $x$. 

Since $[G,L]\leqslant N$, the quotient $K/H$ is cyclic and
generated by $xH$. Let $p\colon K\to K/H$ denote the quotient map. The map $p$ virtually splits. Indeed, if $K/H$ is infinite, the homomorphism $n\mapsto x^n$ defines a section after identifying $K/H$ with $\ZZ$; if $K/H$ is finite, the trivial subgroup of $K/H$ has finite index and admits a section. Since $\psi\in\QQQ(H)^G$ is, in particular, $K$-invariant, the extension theorem \cite[Proposition~6.4]{mixed_Bavard} yields $\nu\in\QQQ(K)$ such that $\nu|_H=\psi$.

Now we show that $\nu \colon K \to \RR$ is $G$-invariant. It follows from $[G,L] \leqslant N$ and the assumption that $ \DDD_{G,K}(\nu)\le  \DDD_{G,N,L}(\mu) < \infty$. Hence Lemma~\ref{lemma ch inv} implies that $\nu$ is $G$-invariant. Hence $(\nu, K) \in \FF$ and $(\nu, K)$ is strictly larger than $(\psi, H)$. This contradicts that $(\psi, H)$ is maximal in $\FF$.
\end{proof}

\subsection{The core and the proof of Theorem~\ref{theorem B}} \label{subsection core}

In this subsection, we introduce the notion of cores of $G$-invariant quasimorphisms on $ \gamma_{q+1} (G)$, and prove Theorem~\ref{theorem B}. We study when $\mu \in \QQQ(N)^G$ is extendable to $L$ as a $G$-invariant homogeneous quasimorphism under the following setting:

\begin{setting} \label{setting}
Let $q$ be a positive integer, and set $N =  \gamma_{q+1} (G)$ and $L =  \gamma_{q} (G)$. Let $\pi \colon F \to G$ be a surjective group homomorphism from a free group $F$ with kernel $R$. Set $M =  \gamma_{q+1} (F)$ and $K =  \gamma_{q} (F)$. Namely, we have the following diagram:
\begin{equation}\label{diagram=core}
\vcenter{\xymatrix@C=36pt{
& F \ar@{->>}[d]_{\pi} \ar@{}[r]|*{\geqslant} & K =  \gamma_{q} (F) \ar@{}[r]|*{\geqslant} \ar@{->>}[d]_{\pi} & M = \gamma_{q+1} (F) \ar@{->>}[d]_{\pi} \\
F/R \ar@{}[r]|*{\cong} & G \ar@{}[r]|*{\geqslant} & L =  \gamma_{q} (G) \ar@{}[r]|*{\geqslant} & N =  \gamma_{q+1} (G).
}}
\end{equation}
\end{setting}

Assume Setting~\ref{setting}. Consider the pullback $\pi^* \mu = \mu \circ (\pi|_{M}) \in \QQQ(M)^F$. Since $F/M = F /  \gamma_{q+1} (F)$ is nilpotent, Corollary~\ref{corollary decomposition} gives $h \in \HHH^1(M)^F$ and $\phi \in \QQQ(F)$ with $\pi^*\mu = h + \phi|_M$. Setting $\nu = \phi|_K \in \QQQ(K)^F$, we have
\[ \pi^* \mu = h + \nu |_M.\]
We will focus on this $F$-invariant homomorphism $h$.

We now show that $h$ does not depend on the decomposition $\pi^* \mu = h + \nu|_M$. Indeed, suppose that $\pi^* \mu = h' + \nu' |_M$ is another decomposition with $h' \in \HHH^1(M)^F$ and $\nu' \in \QQQ(K)^F$. Then we have $h - h' = (\nu' - \nu)|_M$. Since $[F,K] = M$ under Setting~\ref{setting}, if $h \ne h'$, there exist $\tg \in F$ and $\tx \in K$ such that $(h - h')([\tg, \tx]) \ne 0$. Since $h - h'$ is an $F$-invariant homomorphism on $M$, we have
\[ |(\nu' - \nu)([\tg, \tx^n])| = |(h - h')([\tg, \tx^n])| = |n| \cdot |(h - h')([\tg, \tx])|\]
for every $n \in \ZZ$. However, the left-hand side of the above equation is bounded since $\nu' - \nu$ is an $F$-invariant quasimorphism on $K$. This is a contradiction.
 Thus,
\[
\tsigma_1\colon\QQQ(N)^G\to\HHH^1(M)^F,
\quad \mu\mapsto h,
\]
is a well-defined $\RR$-linear map.  In the following proposition, $i_{N,L}$ denotes the inclusion $N \hookrightarrow L$ and $i_{N,L}^*$ is the map $\QQQ(L)^G \to \QQQ(N)^G$, $\phi \mapsto \phi|_N$.

\begin{proposition} \label{proposition tsigma_1}
Assume Setting~\ref{setting}. Then, we have $\Ker (\tsigma_1) = i_{N,L}^* \QQQ(L)^G$.
\end{proposition}
\begin{proof}
Suppose $\mu \in i_{N,L}^* \QQQ(L)^G$. Let $\phi \in \QQQ(L)^G$ such that $\phi|_N = \mu$. Set $\pi^* \phi = \phi \circ (\pi|_K) \in \QQQ(K)^F$. Then $\pi^* \mu = \pi^* \phi |_M$. This implies that $\tsigma_1(\mu) = 0$. Thus we have $i^*_{N,L} \QQQ(L)^G \subset \Ker(\tsigma_1)$.

Next suppose $\mu \in \Ker(\tsigma_1)$. This means that there exists $\nu \in \QQQ(K)^F$ such that $\nu|_M = \pi^* \mu$. Let $g \in G$ and $x \in L$. Let $\tg \in F$ and $\tx \in K$ such that $\pi(\tg) = g$ and $\pi(\tx) = x$. Then we have
\[ |\mu([g,x])| = |\pi^* \mu ([\tg, \tx])| = |\nu([\tg,\tx])| \le \DDD(\nu) < \infty \]
since $\nu \in \QQQ(K)^F$. Then, Proposition~\ref{proposition ch ext} implies that $\mu \in i_{N,L}^* \QQQ(L)^G$. This completes the proof.
\end{proof}

Therefore, the map $\tsigma_1 \colon \QQQ(N)^G \to \HHH^1(M)^F$ induces an injective $\RR$-linear map
\[ \QQQ(N)^G / i^*_{N,L}\QQQ(L)^G \to \HHH^1(M)^F.\]
To obtain an injective $\RR$-linear map from $\WW(G,L,N) = \QQQ(N)^G / (\HHH^1(N)^G + i^*_{N,L} \QQQ(L)^G)$, we need to consider the restriction of $\tsigma_1(\mu)$ to $R \cap M$. Here, recall that $R$ denotes the kernel of $\pi \colon F \to G$ (see Setting~\ref{setting}). Define $\tsigma \colon \QQQ(N)^G \to \HHH^1(R \cap M)^F$ to be the composition
\[ \QQQ(N)^G \xrightarrow{\tsigma_1} \HHH^1(M)^F \xrightarrow{i_{R \cap M,  M }^*} \HHH^1(R \cap M)^F, \]
where $i^*_{R \cap M, M} \colon \HHH^1(M)^{ F } \to \HHH^1(R \cap M)^{ F }$ is the restriction $h \mapsto h|_{R \cap M}$.

\begin{definition} \label{definition core}
Assume Setting~\ref{setting}. Let $\mu \in \QQQ(N)^G$. The \emph{core} of $\mu$ is $\tsigma_1 (\mu)|_{R \cap M} \in \HHH^1(R \cap M)^F$. The map $\tsigma \colon \QQQ(N)^G \to \HHH^1(R \cap M)^F$ is called the \emph{core extractor}.
\end{definition}

\begin{theorem} \label{theorem core kernel}
Assume Setting~\ref{setting}. Then we have $\Ker (\tsigma) = \HHH^1(N)^G + i^*_{N,L} \QQQ(L)^G$. Therefore, if $\mu$ is non-zero in $\WW(G,L,N)$, then there exists $\tz \in R \cap M$ such that $\tsigma(\mu)(\tz) \ne 0$. Namely, the core extractor $\tsigma \colon \QQQ(N)^G \to \HHH^1(R \cap M)^F$ induces an injective $\RR$-linear map $\WW(G,L,N) \to \HHH^1(R \cap M)^F$.
\end{theorem}
\begin{proof}
It is clear that $\HHH^1(N)^G \subset \Ker (\tsigma)$. Proposition~\ref{proposition tsigma_1} implies $i^*_{N,L} \QQQ(L)^G = \Ker (\tsigma_1) \subset \Ker(\tsigma)$. This yields $\HHH^1(N)^G + i_{N,L}^* \QQQ(L)^G \subset \Ker(\tsigma)$.

Now we show $\Ker (\tsigma) \subset \HHH^1(N)^G + i^*_{N,L} \QQQ(L)^G$. Let $\mu \in \Ker (\tsigma)$. Let $h \in \HHH^1(M)^F$ and $\nu \in \QQQ(K)^F$ such that $\pi^*\mu = h + \nu|_{M}$. Since $\mu \in \Ker(\tsigma)$, we have $h|_{R \cap M} = 0$. Then, there exists $\bar{h} \in \HHH^1(N)^G$ such that $h = \bar{h} \circ \pi|_M$. Then we have $(\mu - \bar{h}) \circ \pi|_M = \nu|_M$. Let $g \in G$ and $x \in L$. Let $\tg \in F$, $\tx \in K$ such that $\pi(\tg) = g$ and $\pi(\tx) = x$. Then we have
\[ |(\mu - \bar{h})([g,x])| = |(\mu - \bar{h}) \circ \pi ([\tg, \tx])| =|\nu([\tg, \tx])| \le D(\nu).\]
Hence Proposition~\ref{proposition ch ext} implies that $\mu - \bar{h}$ is extendable to $L$ as a $G$-invariant homogeneous quasimorphism. This means $\mu \in \HHH^1(N)^G + i^*_{N,L} \QQQ(L)^G$.
\end{proof}

The homomorphism $\tsigma_1(\mu)$ and the core $\tsigma(\mu)$ depend on the choice of the free presentation $\pi \colon F \to G$. However, by Proposition~\ref{proposition tsigma_1} and Theorem~\ref{theorem core kernel}, the conditions $\tsigma_1(\mu) = 0$ and $\tsigma(\mu) = 0$ are equivalent to $\mu \in i^*_{N,L}\QQQ(L)^G$ and $\mu \in \HHH^1(N)^G + i^*_{N,L}\QQQ(L)^G$, respectively; in particular, they do not depend on this choice.

Now we are ready to complete the proof of Theorem~\ref{theorem B}.

\begin{proof}[Proof of Theorem~\ref{theorem B}]
Item (1) is a special case of Proposition~\ref{proposition tsigma_1} for $q = 1$, and item (2) is a special case of Theorem~\ref{theorem core kernel} for $q = 1$.
\end{proof}

For $\mu \in \QQQ(N)^G$, the value of the core $\tsigma(\mu)$ can be expressed using a limit as follows:

\begin{proposition} \label{proposition limit formula}
Assume Setting~\ref{setting}. Let $\mu \in \QQQ(N)^G$ and $\tz \in R \cap M$. Let $\tg_1, \cdots, \tg_k \in F$ and $\tx_1, \cdots, \tx_k \in K$ such that $\tz = [\tg_1, \tx_1] \cdots [\tg_k, \tx_k]$. Set $g_i = \pi(\tg_i)$ and $x_i = \pi(\tx_i)$ for $i = 1, \cdots, k$. Then we have
\[ \tsigma(\mu)(\tz) = \lim_{n \to \infty} \frac{\mu([g_1, x_1^n] \cdots [g_k, x_k^n])}{n}.\]
\end{proposition}
\begin{proof}
Set $\tz^{(n)} = [\tg_1, \tx_1^n] \cdots [\tg_k, \tx_k^n]$. Let $h \in \HHH^1(M)^F$ and $\nu \in \QQQ(K)^F$ such that $\pi^* \mu = h + \nu|_M$. Since
\[ [\tg_i, \tx_i^n] = [\tg_i, \tx_i] \cdot \tx_i [\tg_i, \tx_i] \tx_i^{-1} \cdots \tx_i^{n-1}[\tg_i, \tx_i] \tx_i^{-(n-1)}, \]
we have
\[ h(\tz^{(n)}) = h([\tg_1, \tx_1^n] \cdots [\tg_k, \tx_k^n]) = n(h([\tg_1, \tx_1]) + \cdots + h([\tg_k, \tx_k])) = n \cdot h(\tz).\]
Since $\tz^{(n)}$ is a product of $k$ $(F,K)$-commutators of $F$ and $\nu \in \QQQ(K)^F$, we have
\[ |\nu(\tz^{(n)})| \le (2k-1) D(\nu).\]
Hence we have
\[ \lim_{n \to \infty} \frac{\mu([g_1, x_1^n] \cdots [g_k, x_k^n])}{n} = \lim_{n \to \infty} \frac{h(\tz^{(n)}) + \nu(\tz^{(n)})}{n} = h(\tz) = \tsigma(\mu)(\tz). \qedhere\]
\end{proof}

\subsection{Proof of Theorem~\ref{theorem C}} \label{subsection 5.3}

In this subsection, we assume Setting~\ref{setting} and that $\dim_\RR \WW(G,L,N) = n < \infty$. Under these assumptions, we will show in Corollary~\ref{corollary basis 1} that $\Ker(\bT_{N,L})$ admits a basis consisting of elements as described in Proposition~\ref{proposition kernel element}. By combining this with the results from the previous section, we then prove that the coarse kernel of $\iota^G_{N,L} \colon ([G,N], \dGN) \to ([G,L], \dGL)$ is coarsely isomorphic to $(\ZZ^n, \| \cdot \|_1)$ (see Theorem~\ref{theorem basis 2} below). By a direct argument, we can see that this coarse kernel is a splitting coarse kernel (Proposition~\ref{proposition splitting}). In this subsection, we complete the proof of Theorem~\ref{theorem C}.

\begin{theorem} \label{theorem overflow}
Assume Setting~\ref{setting}. Let $\mu \in \QQQ(N)^G$ such that $\mu$ is non-zero in $\WW(G,L,N)$. Then there exist $g_1, \cdots, g_k \in G$ and $x_1, \cdots, x_k \in L$ which satisfy the following properties:
\begin{enumerate}[$(1)$]
\item $[g_1, x_1] \cdots [g_k, x_k] = e_G$.

\item Define $\alpha \in \bCh(G,N)$ by
\[ \alpha = \nlim{\bscl_{G,N}}_{m \to \infty} \left( \frac{[g_1, x_1^m]}{m} + \cdots + \frac{[g_k, x_k^m]}{m} \right).\]
Then, $\mu(\alpha) \ne 0$.
\end{enumerate}
\end{theorem}

\begin{proof}
Let $h \in \HHH^1(M)^F$ and $\nu \in \QQQ( K )^F$ so that $\pi^* \mu = h + \nu|_M$. By Theorem~\ref{theorem core kernel}, there exists $\tz \in R \cap M$ such that $h(\tz) = \tsigma(\mu)(\tz) \ne 0$. Since $M =  \gamma_{q+1} (F)$ and $K =  \gamma_{q} (F)$, there exist $\tg_1, \cdots, \tg_k \in F$ and $\tx_1, \cdots, \tx_k \in K$ such that $\tz = [\tg_1, \tx_1] \cdots [\tg_k, \tx_k]$. Set $g_i = \pi(\tg_i)$ and $x_i = \pi(\tx_i)$ for $i = 1, \cdots, k$.

Then (1) holds since $[g_1, x_1] \cdots [g_k, x_k] = \pi(\tz)$ and $\tz \in R \cap M \subset R = \Ker(\pi)$. 
 For every $m\in\NN$, we have
\[
\left|
 \sum_{i=1}^k\mu([g_i,x_i^m])
 -\mu([g_1,x_1^m]\cdots[g_k,x_k^m])
\right|
\leq (k-1)D(\mu).
\]
Therefore, by Proposition~\ref{proposition limit formula}, condition (2) follows from
\begin{align*}
\mu(\alpha)
&=\lim_{m\to\infty}
  \frac{\sum_{i=1}^k\mu([g_i,x_i^m])}{m}\\
&=\lim_{m\to\infty}
  \frac{\mu([g_1,x_1^m]\cdots[g_k,x_k^m])}{m}\\
&=\tsigma(\mu)(\tz)
 =h(\tz)\ne0. \qedhere
\end{align*}
\end{proof}

\begin{corollary} \label{corollary basis 1}
Assume Setting~\ref{setting} and that $\dim_\RR \WW(G, L ,N) = n$ is finite. Then, there exists a basis $\alpha_1, \dots, \alpha_n$ of $\Ker (\bT_{N, L})$ such that each $\alpha_i$ satisfies the following property $(*)$:
\begin{itemize}
\item[$(*)$] For each $i = 1, \dots, n$, there exist $g_{i,1}, \dots, g_{i, k_i} \in G$ and $x_{i, 1}, \dots, x_{i,k_i} \in L$ such that
\[ [g_{i,1}, x_{i,1}] \cdots [g_{i, k_i}, x_{i, k_i}] = e_G\]
and
\[ \alpha_i = \nlim{\bscl_{G,N}}_{m \to \infty} \left( \frac{[g_{i,1}, x^m_{i,1}]}{m} + \dots + \frac{[g_{i, k_i}, x^m_{i, k_i}]}{m} \right). \]
\end{itemize}
\end{corollary}
\begin{proof}
Let $l$ be an integer such that $0 \le l < n$, and suppose that we have obtained linearly independent elements $\alpha_1, \dots, \alpha_l$ in $\Ker (\bT_{N, L})$ such that each $\alpha_i$ satisfies condition $(*)$. Then, it suffices to show that we can choose an element $\alpha_{l+1}$ satisfying $(*)$ such that $\alpha_1, \dots, \alpha_{l+1}$ are linearly independent. Since $l < n$ and $\Ker(\bT_{N,L})$ is finite-dimensional, there exists $\mu \in \QQQ(N)^G$ such that $\mu$ is non-zero in $\WW(G, L, N)$ and $\mu(\alpha_1) = \dots = \mu(\alpha_l) = 0$ (Theorem~\ref{theorem W duality}). Therefore, by Theorem~\ref{theorem overflow}, there exists an element $\alpha_{l+1}$ satisfying $(*)$ such that $\mu(\alpha_{l + 1}) \ne 0$. Then $\alpha_1, \dots, \alpha_{l+1}$ are linearly independent, completing the proof.
\end{proof}

With the above preparations, we are now ready to determine the coarse kernel of $\iota^G_{N,L} \colon ([G,N], \dGN) \to ([G,L], \dGL)$. Here, we will see that $\iota^G_{N,L}$  admits a splitting coarse kernel in the following sense. 
\begin{definition} \label{definition splitting}
Let $f \colon G \to H$ be a coarse homomorphism between coarse groups. The coarse kernel $j \colon K \to G$ is said to be \emph{splitting} if there exists a coarse homomorphism $r \colon G \to K$ such that $r \circ  j  \bbsim {\rm id}_K$.
\end{definition}

\begin{theorem} \label{theorem basis 2}
Assume Setting~\textup{\ref{setting}} and that $\dim_\RR \WW(G,L,N) = n$ is finite. Then, the natural map $\iota^G_{N,L} \colon ([G,N], \dGN) \to ([G,L], \dGL)$ has a splitting coarse kernel, which is coarsely isomorphic to $(\ZZ^n, \| \cdot \|_1)$.
\end{theorem}
\begin{proof}
By Corollary~\ref{corollary finite-dimensional}, $\Imm(\bT_{N,L})$ is closed.
Let $\alpha_1, \cdots, \alpha_n$ be the basis in $\Ker (\bT_{N,L})$ described in Corollary~\ref{corollary basis 1}. Define the coarse homomorphism $\Phi \colon (\ZZ^n, \| \cdot \|_1) \to ([G,N], \dGN)$ to be the map
\[ \Phi(m_1, \cdots, m_n) = \alpha_1^{(m_1)} \cdots \alpha_n^{(m_n)}\]
described in Proposition~\ref{proposition group lifting}. Let $\Phi'\colon\RR^n\to\bCh(G,N)$ be the linear map defined by
\[
  \Phi'(t_1, \dots, t_n) = t_1 \alpha_1 + \dots + t_n \alpha_n.
\]
Since $\alpha_1,\ldots,\alpha_n$ form a basis of $\Ker(\bT_{N,L})$,
$\Phi'$ is a linear isomorphism onto this kernel.
By Proposition~\ref{proposition group lifting}, $\emb\circ\Phi$
is close to the restriction of $\Phi'$ to $\ZZ^n$.
Here, $\emb \colon [G,N] \to \bCh(G,N)$ is the natural map $x \mapsto x$.

Thus we have the following diagram:
\[ \xymatrix{
(\ZZ^n, \| \cdot \|_1) \ar[r]^-{\Phi} \ar[d] & ([G,N], \dGN) \ar[r]^-{\iota^G_{N,L}} \ar[d]^{\emb} & ([G,L], \dGL) \ar[d]^{\emb} \\
(\RR^n, \| \cdot \|_1) \ar[r]^-{\Phi'}  & (\bCh(G,N), \bscl_{G,N}) \ar[r]^-{\bT_{N,L}} & (\bCh(G,L), \bscl_{G,L}).
}\]
Here, the left vertical arrow is the inclusion. This left square commutes up to closeness by Proposition~\ref{proposition group lifting}, and the right square strictly commutes. Consequently, the proof is completed by Theorem~\ref{theorem Banach} and Propositions~\ref{proposition 1} and \ref{proposition splitting} below.
\end{proof}

\begin{theorem}[{see \cite[Theorem~12.2.5]{LV} and \cite[Theorem 3.15]{fa}}] \label{theorem Banach}
Let $T \colon (X,\| \cdot \|_X) \to (Y, \| \cdot \|_Y)$ be a bounded linear operator between Banach spaces. If $T(X)$ is a closed subspace of $Y$, then  $(\Ker (T), \| \cdot \|_X) $ is the coarse kernel of $T$. 
\end{theorem}

\begin{proposition} \label{proposition 1}
Let
\[ \xymatrix{
& K' \ar[d]^{k} \ar[r]^{u'} & G' \ar[d]^{i} \ar[r]^{f'} & H' \ar[d]^j \\
1 \ar[r] & K \ar[r]^u & G \ar[r]^f & H
}\]
be a diagram of coarse groups that commutes up to closeness. Assume the following:
\begin{enumerate}[$(1)$]
\item The lower horizontal sequence is exact.

\item $i$ and $j$ are coarse embeddings, and $k$ is a coarse equivalence.
\end{enumerate}
Then $u' \colon K' \to G'$ is the coarse kernel of $f'$. Furthermore, if $u \colon K \to G$ is a splitting coarse kernel of $f$, then $u' \colon K' \to G'$ is a splitting coarse kernel of $f' \colon G' \to H'$.
\end{proposition}
\begin{proof}
By the commutativity of the diagram, the composition of $K' \xrightarrow{u'} G' \xrightarrow{f'} H' \xrightarrow{j} H$ is close to the constant map $K' \to H$ at $e_H$. Since $j \colon H' \to H$ is a coarse embedding, the composition of $K' \xrightarrow{u'} G' \xrightarrow{f'} H'$ is close to the constant map $K' \to H'$ at $e_{H'}$. Since $uk \bbsim iu'$ and $uk$ is a coarse embedding, $u'$ is a coarse embedding.

Let $X$ be a coarse space and $\varphi \colon X \to G'$ a bornologous map such that $f' \circ \varphi \bbsim e_{H'}$. Since $jf' \varphi \bbsim fi \varphi \bbsim e_{H}$ and $uk \colon K' \to G$ is the coarse kernel, there exists $\varphi' \colon X \to K'$ such that
\[ i\varphi \bbsim uk \varphi' \bbsim iu' \varphi'.\]
Since $i$ is a coarse embedding, we have $\varphi \bbsim u' \varphi'$. This shows that $u'$ is the coarse kernel of $f'$.

Suppose that $u \colon K \to G$ is a splitting coarse kernel of $f \colon G \to H$. Let $r \colon G \to K$ be a coarse homomorphism satisfying $ru \bbsim {\rm id}_K$. Let $k' \colon K \to K'$ be the coarse inverse of the coarse equivalence $k$. Then $k' \circ r \circ i \colon G' \to K'$ is a coarse homomorphism and $(k' \circ r \circ i) \circ u' \bbsim {\rm id}_{K'}$. This completes the proof.
\end{proof}

\begin{proposition} \label{proposition splitting}
Let $N$ and $L$ be normal subgroups of $G$ such that $N \leqslant L$ and $\dim_\RR \WW(G,L,N) = n < \infty$. Let $\Phi \colon \RR^n \to \Ker(\bT_{N,L})$ be a coarse homomorphism which is coarsely surjective. Then, there exists a bounded linear map $\Psi \colon \bCh(G,N) \to \RR^n$ \textup{(}in particular, a coarse homomorphism\textup{)} such that $\Psi \circ j \circ \Phi \bbsim {\rm id}_{\RR^n}$, where $j \colon \Ker(\bT_{N,L}) \hookrightarrow \bCh(G,N)$ is the inclusion.
\end{proposition}
\begin{proof}
By Theorem~\ref{theorem linear}, there exists a linear map $\Phi' \colon \RR^n \to \Ker(\bT_{N,L})$ such that $\Phi \bbsim \Phi'$.
 Since $\Phi$ is coarsely surjective and $\Phi\bbsim\Phi'$, the linear map $\Phi'$ is surjective. Since $\dim_\RR\Ker(\bT_{N,L})=n$ by Corollary~\ref{corollary finite-dimensional}, $\Phi'$ is a linear isomorphism.  Choose $\mu_1, \dots, \mu_n \in \QQQ(N)^G$ such that their equivalence classes $[\mu_1], \dots, [\mu_n] \in \WW(G,L,N) \cong \Ker(\bT_{N,L})^*$ (Theorem~\ref{theorem W duality}) form the dual basis of $\Phi'(e_1), \dots, \Phi'(e_n)$. 
Define $\boldsymbol{\mu} \colon \bCh(G,N) \to \RR^n$ to be $(\mu_1, \cdots, \mu_n)$, i.e.,
\[ \boldsymbol{\mu}(c) = (\mu_1(c), \dots, \mu_n(c)).\]
Then $\boldsymbol{\mu}$ is a coarse homomorphism satisfying $\boldsymbol{\mu} \circ \Phi' = \mathrm{id}_{\RR^n}$. Therefore, we obtain $\boldsymbol{\mu} \circ \Phi \bbsim \boldsymbol{\mu} \circ \Phi' = \mathrm{id}_{\RR^n}$. Setting $\Psi = \boldsymbol{\mu}$ completes the proof.
\end{proof}

Now we have completed the proof of Theorem~\ref{theorem basis 2}. For the coarse surjectivity in Theorem~\ref{theorem C}, we use the following criterion.

\begin{proposition} \label{proposition xi equivalence}
Assume that $[G, L] / [G,N]$ is a finitely generated abelian group, and that $[L, [G,L]] \leqslant [G,N]$. Then
\[ \iota^G_{N,L} \colon ([G,N], \dGN) \to ([G,L], \dGL)\]
is coarsely surjective.
\end{proposition}
\begin{proof}
Let $u_1, \cdots, u_r$ be elements in $[G,L]$ which form a generating set for $[G,L] / [G,N]$. For each $i = 1, \cdots, r$, let $g_{i,1}, \cdots, g_{i,k_i} \in G$ and $x_{i,1}, \cdots, x_{i,k_i} \in L$ such that
\[ u_i = [g_{i,1}, x_{i,1}] \cdots [g_{i,k_i}, x_{i,k_i}].\]
For $m \in \ZZ$, set
\[ u_i^{(m)} = [g_{i,1}, x_{i,1}^m] \cdots [g_{i,k_i}, x_{i,k_i}^m].\]
Then, in $[G,L]/[G,N]$, we have $u_i^{(m)} = u_i^m$.
Fix $i$ and write $g_j = g_{i,j}$ and $x_j = x_{i,j}$.
Indeed, since $[L,[G,L]] \leqslant [G,N]$, the elements $x_j \in L$ and $[g_j,x_j] \in [G,L]$ commute modulo $[G,N]$.
Since $g_jx_jg_j^{-1}=[g_j,x_j]x_j$, it follows that
\[
g_jx_j^mg_j^{-1}
 =([g_j,x_j]x_j)^m
 =[g_j,x_j]^mx_j^m
 \qquad\text{modulo } [G,N]
\]
for every $m\in\ZZ$. Hence $[g_j,x_j^m]=[g_j,x_j]^m$ in $[G,L]/[G,N]$.

Let $a \in [G,L]$. Then there exist $m_1, \cdots, m_r \in \ZZ$ such that $a (u_1^{(m_1)} \cdots u_r^{(m_r)})^{-1} \in [G,N]$. Hence, there exists $x \in [G,N]$ such that
$ a = x u_1^{(m_1)} \cdots u_r^{(m_r)}.$ Since
\[ d_{\scl_{G,L}}(x,a) \le \cl_{G,L}(x^{-1}a) \le k_1+\cdots+k_r,  \]
$\iota^G_{N,L}$ is coarsely surjective.
\end{proof}

\begin{corollary} \label{corollary xi equivalence}
Suppose that $G$ is finitely generated and that $L =  \gamma_{q} (G)$ and $N =  \gamma_{q+1} (G)$. Then $\iota^G_{N, L} \colon ([G,N], \dGN) \to ([G,L], \dGL)$ is coarsely surjective.
\end{corollary}
\begin{proof}
Since a subgroup of a finitely generated nilpotent group is finitely generated \cite[5.2.18]{Robinson}, $[G,L] / [G,N] = [G,  \gamma_{q} (G)] / [G,  \gamma_{q+1} (G)] =  \gamma_{q+1} (G) /  \gamma_{q+2} (G)$ is finitely generated. Since
\[ [L, [G,L]] = [ \gamma_{q} (G), [G,  \gamma_{q} (G)]] \leqslant [ \gamma_{q} (G),  \gamma_{q+1} (G)] \leqslant [G, \gamma_{q+1} (G)] = [G, N],\]
we have $[L, [G,L]] \leqslant [G,N]$. The assertion follows from Proposition~\ref{proposition xi equivalence}.
\end{proof}

We end this subsection by proving Theorem~\ref{theorem C}.

\begin{proof}[Proof of Theorem~\ref{theorem C}]
By Corollary~\ref{corollary xi equivalence}, $\iota^G_{N,L}$ is coarsely surjective. Since $G$ is finitely generated and $\Gamma = G/N = G/\gamma_{q+1}(G)$ is nilpotent, $\WW(G,L,N)$ is finite-dimensional by Corollary~\ref{corollary fg nilpotent W}~\textup{(2)}. The remaining assertion therefore follows from Theorem~\ref{theorem basis 2}.
\end{proof}

\section{Proof of Theorem~\ref{theorem A}} \label{section main theorem}

In this section, we prove Theorem~\ref{theorem A}. Since $\Gamma=G/N$ is nilpotent, we may choose a positive integer $k$ such that $\gamma_k(G)\leqslant N$. The proof is divided into the following two parts.

\begin{enumerate}[(A)]
\item We combine the splitting coarse kernels of the maps $([G,\gamma_{i+1}(G)], d_{\scl_{G,\gamma_{i+1}(G)}}) \to ([G,\gamma_i(G)], d_{\scl_{G,\gamma_i(G)}})$ to obtain a splitting coarse kernel of $([G,\gamma_k(G)], d_{\scl_{G,\gamma_k(G)}}) \to ([G,G], d_{\scl_G})$ (Subsection~\ref{subsection application}), using general statements on coarse groups proved in Subsection~\ref{subsection general}.

\item We extract the part of this coarse kernel corresponding to $\Ker(\bT_{N,G})$ and realize it as a splitting coarse kernel of $\iota^G_N$ (Subsection~\ref{subsection proof A}).
\end{enumerate}

Subsection~\ref{subsection general} concerns general coarse groups and
requires no assumption on $G$ or $N$. Subsection~\ref{subsection application},
which treats the terms of the lower central series, uses only the finite
generation of $G$; the nilpotency of $\Gamma$ is used only in
Subsection~\ref{subsection proof A}, where we pass from $\gamma_k(G)$ to $N$.

\subsection{Some general statements of coarse groups} \label{subsection general}

In this subsection, we prove several statements about coarse groups which will be used in the next subsection.

\begin{theorem} \label{theorem abelian split}
Let
\[ 0 \to G_0 \xrightarrow{i} G_1 \xrightarrow{p} G_2 \to 0 \]
be a short exact sequence of coarsely abelian groups. Then the following hold:
\begin{enumerate}[$(1)$]
\item Let $r \colon G_1 \to G_0$ be a coarse homomorphism such that $ri \bbsim {\rm id}_{G_0}$. Then there exists a coarse homomorphism $s \colon G_2 \to G_1$ such that $ps \bbsim {\rm id}_{G_2}$ and $ir + sp \bbsim {\rm id}_{G_1}$.

\item Let $s \colon G_2 \to G_1$ be a coarse homomorphism such that $ps \bbsim {\rm id}_{G_2}$. Then there exists a coarse homomorphism $r \colon G_1 \to G_0$ such that $ri \bbsim {\rm id}_{G_0}$ and $ir + sp \bbsim {\rm id}_{G_1}$.

\item If there exist coarse homomorphisms $r \colon G_1 \to G_0$ and $s \colon G_2 \to G_1$ such that $ri \bbsim {\rm id}_{G_0}$, $ps \bbsim {\rm id}_{G_2}$ and $ir + sp \bbsim {\rm id}_{G_1}$, then $rs \bbsim 0_{G_0}$ and the map $(r,p) \colon G_1 \to G_0 \times G_2$ is a coarse equivalence.
\end{enumerate}
\end{theorem}
\begin{proof}
The proof is almost the same as the case of abelian groups.

We now prove (1). Let $r \colon G_1 \to G_0$ be a coarse homomorphism such that $r i \bbsim {\rm id}_{G_0}$. Define $\tilde{s} \colon G_1 \to G_1$ by $\tilde{s} = {\rm id}_{G_1} - ir$. Then, we have
\[ \tilde{s}i = i - iri \bbsim i - i \bbsim 0_{G_{ 1 } }\]
uniformly. Hence we have $\tilde{s}(i(G_0)) \asymp \{0_{G_{ 1 }}\}$. Since $G_1$ is coarsely abelian, $\tilde{s}$ is a coarse homomorphism (see Lemma~\ref{lemma homomorphism sum}). Set $K=i(G_0)$ and let $q\colon G_1\to G_1/K$ be the coarse
quotient map. A coarse homomorphism out of $G_1$ factors through $q$,
up to closeness, if it sends $K$ to a bounded set
\cite[Subsection~5.4]{LV}. Since $\tilde{s}i\bbsim 0_{G_1}$,
this applies to $\tilde{s}$.
Moreover, $K$ realizes the coarse kernel of $p$, and $p$ is
coarsely surjective. Hence the coarse first isomorphism theorem
\cite[Subsection~7.2]{LV} gives a coarse isomorphism
$\bar p\colon G_1/K\to G_2$ such that $\bar p q\bbsim p$.
Via this isomorphism, the factorization of $\tilde{s}$ yields
a coarse homomorphism $s\colon G_2\to G_1$ with
$s p\bbsim\tilde{s}$.
Since $p$ is coarsely surjective, there exist $E \in \EE_{G_2}$ and $\varphi \colon G_2 \to G_1$ such that $(z, p(\varphi(z))) \in E$ for every $z \in G_2$. Then we have
\[ p s \bbsim p s p \varphi \bbsim p \tilde{s} \varphi \bbsim p(\varphi - ir \varphi) \bbsim p \varphi \bbsim {\rm id}_{G_2}.\]
This means $p s \bbsim {\rm id}_{G_2}$. Since
\[ sp \bbsim \tilde{s} = {\rm id}_{G_1} - ir,\]
we have $ir + sp \bbsim {\rm id}_{G_1}$. This completes the proof of (1).

We next prove (2). Suppose that there exists $s \colon G_2 \to G_1$ such that $ps \bbsim {\rm id}_{G_2}$.
Define the function $\tilde{r} \colon G_1 \to G_1$ by $\tilde{r} = {\rm id}_{G_1} - s p$. Since $G_1$ is a coarsely abelian group, $\tilde{r}$ is a coarse homomorphism and 
\[ p \tilde{r} = p ({\rm id}_{G_1} - s p) \bbsim p - p s p \bbsim p -  p \bbsim 0_{G_{ 2 }}.\]
Since $i \colon G_0 \to G_1$ is the coarse kernel of $p \colon G_1 \to G_2$, there exists $r \colon G_1 \to G_0$ such that $ir \bbsim \tilde{r}$. Hence, we have
\[ i r i \bbsim \tilde{r} i = {\rm id}_{G_1} \circ i  - s p i \bbsim i.\]
Since $i$ is a coarse embedding, we have $r i \bbsim {\rm id}_{G_0}$. Since $ir \bbsim \tilde{r} \bbsim {\rm id}_{G_1} - sp$, we have $ir + sp  \bbsim  {\rm id}_{G_1}$. This completes the proof of (2).

Finally, we show (3). First, we show $rs \bbsim 0_{G_0}$ under this condition. Since $ir + sp \bbsim {\rm id}_{G_1}$, we have
\[ r \bbsim rir + rsp \bbsim r + rsp. \]
Since $p$ is coarsely surjective, $rs \bbsim 0_{G_0}$. Next, define $\Phi \colon G_1 \to G_0 \times G_2$ by $\Phi = (r,p)$, and $\Psi \colon G_0 \times G_2 \to G_1$ by $\Psi(x,y) = i(x) + s(y)$. Then, we have the following closeness relations:
\[ \Psi \circ \Phi = ir + sp \bbsim {\rm id}_{G_1}, \]
and
\[ \Phi \circ \Psi(x,y) = \Phi(i(x) + s(y)) \bbsim (ri(x) + rs(y), pi(x) + ps(y)) \bbsim (x,y)\quad ((x,y) \in G_0 \times G_2).\]
Thus $\Phi = (r,p)$ is a coarse isomorphism.
\end{proof}

The next goal is Proposition~\ref{proposition exact sequence} below. 
We use the following two lemmas for the proof of it.

\begin{lemma} \label{lemma a1}
Let $K \xrightarrow{i} G_0 \xrightarrow{f_0} G_1 \xrightarrow{f_1} G_2$ be a sequence of coarse groups, and assume that $f_1$ is a coarse embedding. Then the following are equivalent:
\begin{enumerate}[$(1)$]
\item $i \colon K \to G_0$ is the coarse kernel of $f_0 \colon G_0 \to G_1$.

\item $i \colon K \to G_0$ is the coarse kernel of $f_1 \circ f_0 \colon G_0 \to G_2$.
\end{enumerate}
\end{lemma}
\begin{proof}
In this proof, we can assume that $K$ is a coarse subgroup of $G_0$ and $i \colon K \to G_0$ is the inclusion.

Suppose (1). Let $A$ be a subset of $G_0$ such that $f_1 \circ f_0(A) \preccurlyeq \{ e_{G_2}\}$. Since $f_1$ is a coarse embedding, we have $f_0(A) \preccurlyeq \{ e_{G_1}\}$. Since $K$ is the coarse kernel of $f_0$, we have $A \preccurlyeq K$. Hence, $K$ is the coarse kernel of $f_1 \circ f_0$.

Next, suppose (2). Since $f_1 \circ f_0(K) \preccurlyeq \{ e_{G_2}\}$ and $f_1$ is a coarse embedding, we have $f_0(K) \preccurlyeq \{ e_{G_1}\}$. Let $A$ be a subset of $G_0$ such that $f_0(A) \preccurlyeq \{ e_{G_1}\}$. As $f_1 \circ f_0(A) \preccurlyeq \{ e_{G_2}\}$ and $K$ is the coarse kernel of $f_1 \circ f_0$, we have $A \preccurlyeq K$. Thus $K$ is the coarse kernel of $f_0$. This completes the proof.
\end{proof}

\begin{lemma} \label{lemma a2}
Let $G_0 \xrightarrow{f_0} G_1 \xrightarrow{f_1} G_2$ be a sequence of coarse groups such that $f_0$ is a coarse embedding. Let $i \colon K \to G_1$ be the coarse kernel of $f_1$, and $j \colon K \to G_0$ a coarse homomorphism such that $f_0 \circ j \bbsim i$. Then $j \colon K \to G_0$ is the coarse kernel of $f_1 \circ f_0$.
\end{lemma}
\begin{proof}
Since $i \bbsim f_0 \circ j$ and $i$ is a coarse embedding, $j$ is a coarse embedding. Moreover, $f_1 \circ f_0 \circ j \bbsim f_1 \circ i \bbsim e_{G_2}$. Let $X$ be a coarse space and let $\varphi \colon X \to G_0$ be a bornologous map such that $f_1 \circ f_0 (\varphi(X)) \preccurlyeq \{ e_{G_2}\}$. Since $i \colon K \to G_1$ is the coarse kernel of $f_1$, there exists $\varphi' \colon X \to K$ satisfying $i \circ \varphi' \bbsim f_0 \circ \varphi$. Hence we have
\[ f_0 \circ \varphi \bbsim i \circ \varphi' \bbsim f_0 \circ j \circ \varphi'.\]
Since $f_0$ is a coarse embedding, we have $\varphi \bbsim j \circ \varphi'$. This implies that $j \colon K \to G_0$ is the coarse kernel of $f_1 \circ f_0$.
\end{proof}

\begin{proposition} \label{proposition exact sequence}
Let $G_0 \xrightarrow{f_0} G_1 \xrightarrow{f_1} G_2$ be a sequence of coarse groups, and suppose that $f_0$ is coarsely surjective. If $f_0$ and $f_1$ have coarse kernels $\KKer(f_0)$ and $\KKer(f_1)$, respectively, then $f_1 \circ f_0$ also has a coarse kernel $\KKer(f_1 \circ f_0)$, and the following sequence
\begin{align} \label{shortex}
1 \to \KKer(f_0) \to \KKer(f_1 \circ f_0) \xrightarrow{\check{f}_0} \KKer(f_1) \to 1
\end{align}
is an exact sequence. Here, each arrow in this short exact sequence is induced by the universality of the coarse kernels. Furthermore, if $G_0$ is coarsely abelian and the coarse kernels of $f_0$ and $f_1$ are splitting coarse kernels, then the coarse kernel of $f_1 \circ f_0$ is also a splitting coarse kernel.
\end{proposition}
\begin{proof}
Let $K_1 \subset G_1$ be the realization of $\KKer(f_1)$, i.e., $K_1$ is a subset of $G_1$ such that the inclusion $K_1 \hookrightarrow G_1$ is the coarse kernel. Similarly, let $K_0 \subset G_0$ be a realization of $\KKer(f_0)$. We can assume that $e_{G_1} \in K_1$ and $e_{G_0} \in K_0$. Since $f_0$ is coarsely surjective, there exists $E \in \EE_{G_1}$ satisfying the following conditions:
\begin{itemize}
\item $E$ is a symmetric subset of $G_1 \times G_1$ containing the diagonal $\Delta_{G_1}$.

\item $E[f_0(G_0)] = G_1$.

\item $f_0( K_0 ) \subset E_{e_{G_1}}$.

\item $(f_0(e_{G_0}), e_{G_1}) \in E$.
\end{itemize}

Set $K = K_0 \cup \bigl(K_0 \cdot f_0^{-1}(E[K_1])\bigr)$. Then we show that $K$ is the coarse kernel of $f_1 \circ f_0$.

We first show that $f_1 \circ f_0(K) \asymp e_{G_2}$. Since $K_0$ is the coarse kernel of $f_0$, $f_0(K_0) \preccurlyeq \{ e_{G_1}\}$, and hence $f_1 \circ f_0(K_0) \preccurlyeq \{ e_{G_2} \}$. For every $y \in f_0^{-1}(E[K_1])$, there exists $\varphi(y) \in K_1$ such that $f_0(y) \bsim_E\, \varphi(y)$. Since $K_1$ is the coarse kernel of $f_1$ and the inclusion $K_1 \hookrightarrow E[K_1]$ is a coarse equivalence , we have $f_1 \circ f_0( f_0^{-1}(E[K_1])) \preccurlyeq \{e_{G_2}\}$. Since $f_0$ and $f_1$ are coarse homomorphisms, we have the following closeness among functions on $K_0 \times f_0^{-1}(E[K_1])$:
\[ f_1 \circ f_0(xy) \sim f_1 \circ f_0(x) \cdot f_1 \circ f_0(y) \sim e_{G_2} \cdot e_{G_2} \sim e_{G_2}.\]
Together with $f_1 \circ f_0(K_0) \preccurlyeq \{e_{G_2}\}$, this implies $f_1 \circ f_0 (K) \asymp e_{G_2}$.

Next we show that $K$ satisfies the property of the coarse kernel. Let $X$ be a subset of $G_0$ such that $f_1 \circ f_0 (X) \asymp \{ e_{G_2}\}$. Since $K_1$ is the coarse kernel of $f_1$, we have $f_0(X) \preccurlyeq K_1$. Thus there exists $E'\in\EE_{G_1}$ such that $f_0(X)\subset E'[K_1]$.
For each $x\in X$, choose $k_x\in K_1$ such that $(f_0(x),k_x)\in E'$.
Since $E[f_0(G_0)]=G_1$ and $E$ is symmetric, we may choose $\varphi(x)\in G_0$ such that $(f_0(\varphi(x)),k_x)\in E$.
Then $\varphi(x)\in f_0^{-1}(E[K_1])$ and
\[
f_0(\varphi(x))\bsim_{E\circ(E')^{-1}}\,f_0(x)
\qquad (x\in X).
\]
Consequently, the following maps from $X$ to $G_1$ are close:
\[ f_0(x \cdot \varphi(x)^{-1}) \bbsim f_0(x) \cdot f_0(\varphi(x)^{-1}) \bbsim f_0(x) \cdot f_0(\varphi(x))^{-1} \bbsim f_0(x) \cdot f_0(x)^{-1} \bbsim e_{G_1}, \quad (x \in X).\]
Since $K_0$ is the coarse kernel of $f_0$, we have $\{ x \cdot  \varphi(x)^{-1}  \; | \; x \in X\}  \preccurlyeq   K_0 $. Thus there exists a function $\psi \colon X \to K_0$ satisfying the following closeness between the functions on $X$:
\[ x \cdot \varphi(x)^{-1} \bbsim \psi(x).\]
Then we have the following closeness on the functions on $X$:
\[ x \bbsim x \cdot (\varphi(x)^{-1} \cdot \varphi(x)) \bbsim (x \cdot \varphi(x)^{-1}) \cdot \varphi(x) \bbsim \psi(x) \cdot \varphi(x) \in K_0 \cdot f_0^{-1}( E[K_1] ) \subseteq K. \]
This means $X \preccurlyeq  K $. This completes the proof that $K$ is the coarse kernel of $f_1 \circ f_0$.

Next we show the exactness of the sequence \eqref{shortex}. Consider the diagram
\[ \xymatrix{
K_0 \ar[r] \ar[d] & K \ar[r]^{\check{f}_0} \ar[d] & K_1 \ar[d] \\
G_0 \ar[r]^{\rm id} & G_0 \ar[r]^{f_0} & G_1
}\]
which commutes up to closeness. Then, $K \to G_0$ is a coarse embedding and $K_0 \to G_0$ is the coarse kernel of $f_0 \colon G_0 \to G_1$. Hence, by Lemma~\ref{lemma a2}, the map $K_0 \hookrightarrow K$ is the coarse kernel of the composition $K \to G_0 \to G_1$. Since $K_1 \to G_1$ is a coarse embedding, Lemma~\ref{lemma a1} implies that the inclusion $K_0 \hookrightarrow K$ is the coarse kernel of $\check{f}_0 \colon K \to K_1$. Thus the sequence
\[ 1 \to \KKer(f_0) \to \KKer(f_1 \circ f_0) \xrightarrow{\check{f}_0} \KKer(f_1) \]
is exact. By the construction of $K$, it is clear that the map $\check{f}_0 \colon K \to K_1$ is coarsely surjective. This completes the proof that the sequence \eqref{shortex} is exact.

Finally, suppose that $G_0$ is coarsely abelian. Set $\KKer(f_0) = K_0$, $\KKer(f_1 \circ f_0) = K$, and $\KKer(f_1) = K_1$. Consider the following diagram:
\[ \xymatrix{
& G_0 \ar[r]^{{\rm id}_{G_0}} & G_0 \ar[r]^{f_0} & G_1 & \\
0 \ar[r] & K_0 \ar[r]^j \ar[u]_{i_0} & K \ar[r]^p \ar[u]_i &  K_1 \ar[r] \ar[u]_{i_1} & 0.
}\]
Since $G_0$ is coarsely abelian, $K_0$ and $K$ are coarsely abelian. Since $p \colon K \to K_1$ is a coarsely surjective coarse homomorphism from a coarsely abelian group, $K_1$ is also coarsely abelian.

Let $r_k \colon G_k \to K_k$ be a coarse homomorphism such that $r_k i_k \bbsim {\rm id}_{K_k}$ for $k = 0, 1$. Then $r = r_0 i$ is a coarse homomorphism such that $rj = r_0 i j \bbsim r_0 i_0 \bbsim {\rm id}_{K_0}$. Hence Theorem~\ref{theorem abelian split} implies that there exists a coarse homomorphism $s \colon K_1 \to K$ such that $ps \bbsim {\rm id}_{K_1}$ and $ jr  + sp \bbsim {\rm id}_{K}$. Define the coarse homomorphism $G_0 \to K$ by
\[ \Phi = jr_0 + s r_1 f_0.\]
We have to show that $\Phi i \bbsim {\rm id}_K$. To see this, since $(r, p) \colon K \to K_0 \times K_1$ is a coarse equivalence (Theorem~\ref{theorem abelian split}), it suffices to show
\[ (r, p) \circ \Phi \circ i(x) \bbsim (r, p)(x)\]
uniformly with respect to $x \in K$. Indeed, we have
\[ r \circ (j r_0 + s r_1 f_0) \circ i \bbsim r_0 i j r_0 i + r s r_1 f_0 i \bbsim r_0 i_0 r_0 i + 0 \bbsim r_0 i = r,\]
\[ p \circ (j r_0 + s r_1 f_0) \circ i \bbsim 0 + psr_1 f_0 i \bbsim r_1 f_0 i \bbsim r_1 i_1 p \bbsim p.\]
Here we use $rs \bbsim 0_{K_0}$ (see Theorem~\ref{theorem abelian split}). This completes the proof.
\end{proof}

\subsection{Applications to $\scl_{G,N}$} \label{subsection application}

In this subsection, we combine Theorem~\ref{theorem C} with the results in the previous subsection to study the coarse homomorphism $\iota^G_N \colon ([G,N], \dGN) \to ([G,G], \dG)$, when $N = \gamma_q(G)$.

\begin{theorem} \label{theorem direct sum}
Let $G$ be a finitely generated group, $m$ and $n$  positive  integers such that $m \le n$. Set $N = \gamma_n(G)$, $L = \gamma_m(G)$, and set $l = \dim_\RR \WW(G, L, N)$. Then, the coarse homomorphism
\[ \iota^G_{N,L} \colon ([G, N], \dGN) \to ([G,L], \dGL)\]
is coarsely surjective and has a splitting coarse kernel that is coarsely isomorphic to $(\ZZ^l, \| \cdot \|_1)$.
\end{theorem}
\begin{proof}
Since $G$ is finitely generated and $G/N$ is nilpotent, $\WW(G,L,N)$ is finite-dimensional by Corollary~\ref{corollary fg nilpotent W}~\textup{(2)}.
The case $n = m$ is clear. We prove the case $n > m$ by induction with respect to $n - m$. Consider the following sequence:
\[ ([G,\gamma_n(G)], d_{\scl_{G, \gamma_n(G)}}) \xrightarrow{f_0} ([G,\gamma_{n-1}(G)],  d_{\scl_{G, \gamma_{n-1}(G)}}) \xrightarrow{f_1} ([G, \gamma_m(G)], d_{\scl_{G, \gamma_m(G)}}),\]
where $f_0$ and $f_1$ are natural maps. Then $f_1 \circ f_0 = \iota^G_{N,L}$. It follows from Theorem~\ref{theorem C} that $f_0$ is coarsely surjective and has a splitting coarse kernel $K_1$. By the induction hypothesis, $f_1$ is coarsely surjective and has a splitting coarse kernel $K_2$; in particular, $f_1 \circ f_0$ is coarsely surjective. Since these coarse groups are coarsely abelian, it follows from Proposition~\ref{proposition exact sequence} that $f_1 \circ f_0$ has a splitting coarse kernel $K$, and the sequence
\[ 0 \to K_1 \to K \to K_2 \to 0\]
of coarsely abelian groups is exact. Since $K_1 \hookrightarrow [G, \gamma_n(G)]$ is coarsely splitting, we conclude that $K_1 \hookrightarrow K$ is also coarsely splitting. Hence, Theorem~\ref{theorem abelian split} implies that $K$ and $K_1 \times K_2$ are coarsely isomorphic. Since $\WWW(G,\gamma_n(G))$ is finite-dimensional by Proposition~\ref{proposition fg nilpotent}, Proposition~\ref{proposition W injection}~\textup{(2)}, applied to $\gamma_n(G) \leqslant \gamma_{n-1}(G) \leqslant \gamma_m(G)$, gives
\[ \dim_\RR \WW(G, \gamma_m(G), \gamma_n(G)) = \dim_\RR \WW(G,\gamma_m(G), \gamma_{n-1}(G)) + \dim_\RR \WW(G, \gamma_{n-1}(G), \gamma_n(G)).\]
 Together with $K\cong K_1\times K_2$, this dimension formula shows that $K$ is coarsely isomorphic to $(\ZZ^l,\|\cdot\|_1)$. 
\end{proof}

\subsection{Proof of Theorem~\ref{theorem A}} \label{subsection proof A}

The goal of this subsection is to complete the proof of Theorem~\ref{theorem A}.

Let $k$ be a positive integer such that $\gamma_k(G) \leqslant N$. Since $\bT_{N,G} \circ \bT_{\gamma_k(G),N} = \bT_{\gamma_k(G),G}$, the operator $\bT_{\gamma_k(G),N}$ maps $\Ker(\bT_{\gamma_k(G),G})$ into $\Ker(\bT_{N,G})$. Let
\[ \rho \colon \Ker(\bT_{\gamma_k(G),G}) \to \Ker(\bT_{N,G}) \]
denote the restriction of $\bT_{\gamma_k(G),N}$.

\begin{lemma} \label{lemma surjective}
Let $G$ be a finitely generated group such that $ \Gamma = G/N$ is nilpotent. Let $k$ be a  positive  integer such that $\gamma_k(G) \leqslant N$. Then $\rho \colon \Ker(\bT_{\gamma_k(G), G}) \to \Ker(\bT_{N,G})$ is surjective.
\end{lemma}
\begin{proof}
In this case, the spaces $\Ker(\bT_{N,G})$, $\Ker(\bT_{\gamma_k(G), G})$, $\WWW(G,N)$, and $\WWW(G,\gamma_k(G))$ are finite-dimensional (Proposition~\ref{proposition fg nilpotent} and Corollary~\ref{corollary finite-dimensional}).
 By Theorem~\ref{theorem W duality}, we have the natural isomorphisms $\Ker(\bT_{N,G})^* \cong \WWW(G,N)$ and $\Ker(\bT_{\gamma_k(G),G})^* \cong \WWW(G,\gamma_k(G))$. 
Under these isomorphisms, the dual map of $\rho$ is $\WWW(G,N) \to \WWW(G,\gamma_k(G))$, $[\mu] \mapsto [\mu|_{\gamma_k(G)}]$. Hence, it suffices to show that this map is injective, but this follows from Proposition~\ref{proposition W injection}~\textup{(1)}, since $G/N$ and $G/\gamma_k(G)$ are nilpotent, hence amenable and boundedly $3$-acyclic.
\end{proof}

The following comparison lemma identifies an abstract coarse kernel, realized inside a coarse group, with the kernel of a bounded linear operator.

\begin{lemma} \label{lemma factorization}
Let $T \colon X \to Y$ be a bounded linear operator between real Banach spaces such that $T(X)$ is closed, and assume that $\Ker(T)$ has finite dimension $m$. Let $A$ and $K$ be coarse groups such that $K$ is coarsely isomorphic to $(\ZZ^m, \|\cdot\|_1)$, and let $\theta \colon A \to X$ and $j \colon K \to A$ be coarse homomorphisms such that $\theta \circ j$ is a coarse embedding and $T \circ \theta \circ j \bbsim 0_Y$. Then there exists a coarse isomorphism $u \colon K \to \Ker(T)$ such that $\iota \circ u \bbsim \theta \circ j$, where $\iota \colon \Ker(T) \hookrightarrow X$ is the inclusion.
\end{lemma}
\begin{proof}
Composing the coordinatewise floor map $(\RR^m,\|\cdot\|_1) \to (\ZZ^m,\|\cdot\|_1)$ with a coarse isomorphism $(\ZZ^m,\|\cdot\|_1) \to K$, we obtain a coarse isomorphism $\varphi \colon (\RR^m,\|\cdot\|_1) \to K$; let $\psi \colon K \to \RR^m$ be a coarse inverse of $\varphi$, which is again a coarse homomorphism. Set $a = \theta \circ j \circ \varphi$; this is a coarse homomorphism and a coarse embedding with $T \circ a \bbsim 0_Y$. Since $T(X)$ is closed, $\iota$ is the coarse kernel of $T$ by Theorem~\ref{theorem Banach}, so there exists a bornologous map $a' \colon \RR^m \to \Ker(T)$ such that $\iota \circ a' \bbsim a$ (Remark~\ref{remark coarse kernel}). Since $\iota$ is an isometric embedding, $a'$ is a coarse homomorphism and a coarse embedding. Identify $\Ker(T)$ with $(\RR^m, \|\cdot\|_1)$ via a linear isomorphism, which is bi-Lipschitz since $\dim_\RR \Ker(T) = m$. Then Theorem~\ref{theorem linear} gives a linear map $S \colon \RR^m \to \Ker(T)$ with $a' \bbsim S$. Since $S \bbsim a'$, $S$ is a coarse embedding; hence it is injective, and hence a linear isomorphism. In particular, $a'$ is a coarse isomorphism. Therefore $u = a' \circ \psi$ is a coarse isomorphism, and
\[
\iota \circ u = \iota \circ a' \circ \psi \bbsim \theta \circ j \circ \varphi \circ \psi \bbsim \theta \circ j. \qedhere
\]
\end{proof}

Now we are ready to prove Theorem~\ref{theorem A}.

\begin{proof}[Proof of Theorem~\textup{\ref{theorem A}}]
We first prove~$(1)$. By Theorem~\ref{theorem direct sum}, $\iota^G_{\gamma_k(G)}$ is coarsely surjective. Since $\iota^G_{\gamma_k(G)} = \iota^G_N \circ \iota^G_{\gamma_k(G),N}$, so is $\iota^G_N$.

We next prove~$(2)$. Note that $\WWW(G,N)$ and $\WWW(G,\gamma_k(G))$ are finite-dimensional by Proposition~\ref{proposition fg nilpotent}. Set $m = \dim_\RR \WWW(G,\gamma_k(G))$. By Theorem~\ref{theorem direct sum}, $\iota^G_{\gamma_k(G)}$ has a coarse kernel $j \colon K \to [G,\gamma_k(G)]$ such that $K$ is coarsely isomorphic to $(\ZZ^m,\|\cdot\|_1)$. Then $\emb \circ j$ is a coarse embedding, and $\bT_{\gamma_k(G),G} \circ \emb \circ j = \emb \circ \iota^G_{\gamma_k(G)} \circ j \bbsim 0$. By Corollary~\ref{corollary finite-dimensional}, $\Imm(\bT_{\gamma_k(G),G})$ is closed and $\dim_\RR \Ker(\bT_{\gamma_k(G),G}) = m$. Hence Lemma~\ref{lemma factorization}, applied to $T = \bT_{\gamma_k(G),G}$, $A = ([G,\gamma_k(G)], d_{\scl_{G,\gamma_k(G)}})$ and $\theta = \emb$, gives a coarse isomorphism $u \colon K \to \Ker(\bT_{\gamma_k(G),G})$ such that $\iota \circ u \bbsim \emb \circ j$, where $\iota$ is the inclusion. Let $v$ be a coarse inverse of $u$ and set $F = j \circ v \colon \Ker(\bT_{\gamma_k(G),G}) \to [G,\gamma_k(G)]$. Then $F$ is a coarse homomorphism and $\emb \circ F \bbsim \iota \circ u \circ v \bbsim \iota$.

By Lemma~\ref{lemma surjective}, $\rho$ is a linear surjection between finite-dimensional spaces, and hence admits a linear section $\sigma \colon \Ker(\bT_{N,G}) \to \Ker(\bT_{\gamma_k(G),G})$. Set
\[ u' = \iota^G_{\gamma_k(G),N} \circ F \circ \sigma \colon \Ker(\bT_{N,G}) \to ([G,N], \dGN). \]
Since $\emb \circ \iota^G_{\gamma_k(G),N} = \bT_{\gamma_k(G),N} \circ \emb$ and $\bT_{\gamma_k(G),N}$ is Lipschitz, we have
\[ \emb \circ u' = \bT_{\gamma_k(G),N} \circ \emb \circ F \circ \sigma \bbsim \bT_{\gamma_k(G),N} \circ \sigma = \rho \circ \sigma, \]
and the right-hand side is the inclusion $\Ker(\bT_{N,G}) \hookrightarrow \bCh(G,N)$. Thus the diagram
\[ \xymatrix{
 & \Ker(\bT_{N,G}) \ar[d]_{\rm id} \ar[r]^-{u'} & ([G,N], \dGN) \ar[d]^{\emb} \ar[r]^{\iota^G_N} & ([G,G], d_{\scl_G}) \ar[d]^{\emb} \\
0 \ar[r] & \Ker(\bT_{N,G}) \ar[r] & (\bCh(G,N), \bscl_{G,N}) \ar[r]^{\bT_{N,G}} & (\bCh(G,G), \bscl_G)
}\]
commutes up to closeness. The lower horizontal sequence is exact by Corollary~\ref{corollary finite-dimensional} and Theorem~\ref{theorem Banach}, and its coarse kernel is splitting by Proposition~\ref{proposition splitting}, applied with $L = G$ to a linear isomorphism $\Phi \colon \RR^n \to \Ker(\bT_{N,G})$ (the retraction being $\Phi \circ \Psi$). Hence Proposition~\ref{proposition 1} shows that $u'$ is a splitting coarse kernel of $\iota^G_N$. Since $\dim_\RR \Ker(\bT_{N,G}) = n$ by Corollary~\ref{corollary finite-dimensional}, this coarse kernel is coarsely isomorphic to $(\ZZ^n,\|\cdot\|_1)$. This proves~$(2)$.

Finally, identify the splitting coarse kernel in~$(2)$ with $(\ZZ^n,\|\cdot\|_1)$, and let $r$ be the associated retraction. By~$(1)$ and~$(2)$, the sequence
\[ 0 \to (\ZZ^n,\|\cdot\|_1) \to ([G,N],\dGN) \xrightarrow{\iota^G_N} ([G,G],\dG) \to 0 \]
is a short exact sequence of coarsely abelian groups. Hence Theorem~\ref{theorem abelian split}~$(1)$ and~$(3)$ imply that $(r,\iota^G_N)$ is an isomorphism of coarse groups.
\end{proof}

\section{Explicit coarse kernels} \label{section explicit}

In this section, under the assumptions of Theorem~\ref{theorem A}, we give an explicit representation of the coarse kernel of $\iota^G_N$ as a subset of $[G,N]$.

First, as motivation for providing an explicit representation of the coarse kernel,  we discuss the following application of coarse kernels (Proposition~\ref{proposition W vanishing}): given $\mu \in \QQQ(N)^G$, whether it is non-zero in $\WWW(G,N)$ can be determined by whether it is bounded on the coarse kernel.

\begin{proposition} \label{proposition W vanishing}
Let $G$ be a finitely generated group, $N$ a normal subgroup such that $\Gamma=G/N$ is nilpotent, and let $\mu \in \QQQ(N)^G$. Let $K$ be the coarse kernel of $\iota^G_N \colon ([G,N], \dGN) \to ([G,G], \dG)$ realized as a subset of $[G,N]$. Then the following are equivalent:
\begin{enumerate}[$(1)$]
\item $\mu$ vanishes in $\WWW(G,N)$, i.e., $\mu \in \HHH^1(N)^G + i^* \QQQ(G)$.

\item $\mu |_K$ is a bounded function on $K$.
\end{enumerate}
\end{proposition}
\begin{proof}
By Theorem~\ref{theorem A}, $K$ is coarsely isomorphic to $(\ZZ^n,\|\cdot\|_1)$, where $n = \dim_\RR \WWW(G,N)$. Moreover, $\Imm(\bT_{N,G})$ is closed and $\dim_\RR \Ker(\bT_{N,G}) = n$ (Corollary~\ref{corollary finite-dimensional}), and $\bT_{N,G} \circ \emb|_K = \emb \circ \iota^G_N|_K \bbsim 0$. Hence Lemma~\ref{lemma factorization}, applied to $T = \bT_{N,G}$, $A = ([G,N],\dGN)$, $\theta = \emb$ and the inclusion $K \hookrightarrow [G,N]$, gives a coarse isomorphism $u \colon K \to \Ker(\bT_{N,G})$ such that the following diagram commutes up to closeness:
\[ \xymatrix{
0 \ar[r] & K \ar[d]_{u} \ar[r] & ([G,N], \dGN) \ar[d] \ar[r]^{\iota^G_N} & ([G,G], d_{\scl_G}) \ar[d] \\
0 \ar[r] & \Ker(\bT_{N,G}) \ar[r] & (\bCh(G,N), \bscl_{G,N}) \ar[r]^{\bT_{N,G}} & (\bCh(G,G), \bscl_G).
}\]
 Since the left vertical arrow is a coarse equivalence and the diagram commutes up to closeness, $\mu|_K$ is bounded if and only if $\mu|_{\Ker(\bT_{N,G})}$ is bounded. 
Since $\mu|_{\Ker(\bT_{N,G})}$ is linear, it is bounded as a function if and only if it vanishes identically.
By Theorem~\ref{theorem W duality}, the vanishing of $\mu|_{\Ker(\bT_{N,G})}$ is equivalent to the vanishing of the class of $\mu$ in $\WWW(G,N)$.
\end{proof}

\begin{remark}
In our Part I paper \cite[Theorem~2.3]{fa}, we proved such a characterization at the level of rational chains (in \cite[Theorem~B]{fa}, we also showed that the coarse kernel at the level of rational chains exists if and only if  $\widehat{i^*}_{N,G}(\QQQ(G)/\HHH^1(G))$ is norm-closed in $(\QQQ(N)^G/\HHH^1(N)^G,2D_N)$).
In \cite[Theorem~5.2]{fa}, we proved the corresponding group-level
characterization in full generality, but the proof required passing to
further positive-integer multiples and then scaling down. By contrast,
under the assumptions of Theorem~\ref{theorem A}, the coarse kernel
constructed above provides uniformly controlled group realizations of
each multiple $m\alpha$ itself, without any additional scaling factor
(see Lemma~\ref{lemma lifting} and Corollary~\ref{corollary basis 1}),
and the characterization in Proposition~\ref{proposition W vanishing} is stated directly
in terms of this coarse kernel.
\end{remark}

Now, let us consider an explicit expression of the coarse kernel. However, finding such an expression appears to be a challenging problem in general. 
Indeed, when determining the coarse kernel in Setting~\ref{setting}, it was crucial to find an element $\tz \in R \cap M$ for which the value of the core does not vanish in Theorem~\ref{theorem core kernel}. Furthermore, the method for finding the core is not constructive, and thus such an element $\tz$ was not obtained in a constructive manner.

Nevertheless, under certain conditions, it is possible to determine the coarse kernel explicitly. 
We now specialize the core extractor introduced in Subsection~\ref{subsection core} to the case $N=[G,G]$. 
Let $\pi \colon F \to G$ be a surjective group homomorphism from a free group $F$, let $R=\Ker(\pi)$, and set $M=[F,F]$. 
 For $\mu\in\QQQ(N)^G$, let $h_\mu\in\HHH^1(M)^F$ be the unique
homomorphism appearing in a decomposition
\[
\pi^*\mu=h_\mu+\phi|_M
\]
with $\phi\in\QQQ(F)$ (Corollary~\ref{corollary decomposition}). Thus
$\tsigma(\mu)=h_\mu|_{R\cap M}$.

\begin{theorem} \label{theorem explicit}
Let $G$ be a finitely generated group, $N$ the commutator subgroup $[G,G]$ of $G$, and $\pi \colon F \to G$ a surjective group homomorphism from a free group $F$. Set $R = \Ker(\pi)$ and $M = [F,F]$. Let $n = \dim_\RR \WWW(G,N)$, which is finite by Proposition~\ref{proposition fg nilpotent}. 
 Suppose that $\tz_1,\dots,\tz_n\in M$ satisfy
$\pi(\tz_i)\in[G,N]$ for every $i$ and the following condition $(*)$:
\begin{itemize}
\item[$(*)$] For every $\mu\in\QQQ(N)^G$,
\[
h_\mu(\tz_1)=\cdots=h_\mu(\tz_n)=0
\quad\Longrightarrow\quad
\tsigma(\mu)=0.
\]
\end{itemize} 
Choose elements $\ta_{i,1}, \dots, \ta_{i,k_i}, \tb_{i,1}, \dots, \tb_{i,k_i} \in F$ such that
\[ \tz_i = [\ta_{i,1}, \tb_{i,1}] \cdots [\ta_{i, k_i}, \tb_{i, k_i}]. \]
Set $a_{i,j} = \pi(\ta_{i,j})$, $b_{i,j} = \pi(\tb_{i,j})$ and
\[ z_i^{(m)} = [a_{i,1}, b_{i,1}^m] \cdots [a_{i, k_i}, b_{i, k_i}^m]. \]
Then, the set 
\[ \{ z_1^{(m_1)} \cdots z_n^{(m_n)} \mid m_1, \dots, m_n \in \ZZ \} \]
is a realization of the coarse kernel of $\iota^G_N \colon ([G,N], \dGN) \to ([G,G], \dG)$.
\end{theorem}

\begin{remark} \label{remark explicit}
For a subset $S$ of a group $H$, write $\llangle S\rrangle_H$
for the smallest normal subgroup of $H$ containing $S$.
A sufficient condition for $(*)$ is that
$R\cap M=\llangle\tz_1,\ldots,\tz_n\rrangle_F$.
More generally, suppose that
$R\cap M=\llangle r_1,\ldots,r_n\rrangle_F$ and that
$\tz_i$ and $r_i$ have the same image in $F/\gamma_3(F)$
for each $i$.
Then $\pi(\tz_i)\in[G,N]$ and
$h_\mu(\tz_i)=h_\mu(r_i)$, since every $h_\mu\in\HHH^1([F,F])^F$    vanishes on $\gamma_3(F) = [F, [F, F]]$. 
 Thus condition $(*)$ holds. 
\end{remark}

\begin{proof}[Proof of Theorem~\ref{theorem explicit}]
 If $n=0$, the assertion follows immediately from Theorem~\ref{theorem A}. Hence, we assume that $n>0$.  We set
\[ z'_i = \nlim{\bscl_{G,N}}_{m \to \infty} \left( \frac{[a_{i,1}, b_{i,1}^m]}{m} + \cdots + \frac{[a_{i,k_i}, b_{i, k_i}^m]}{m} \right).\]
 Then $z'_i\in\Ker(\bT_{N,G})$ by Proposition~\ref{proposition kernel element}, since $\pi(\tz_i)\in[G,N]$. 
We define $\Phi \colon \ZZ^n \to [G,N]$ by
\[ \Phi(m_1, \dots, m_n) = z_1^{(m_1)} \dots z_n^{(m_n)}, \]
and $\Phi' \colon \ZZ^n \to \Ker(\bT_{N,G})$ by
\[ \Phi'(m_1, \dots, m_n) = m_1 z'_1 + \dots + m_n z'_n.\]
Consider the following diagram:
\begin{align} \label{explicit}
\vcenter{
\xymatrix{
& \ZZ^n \ar[r]^-{\Phi} \ar[d]_{\Phi'} & ([G,N], \dGN) \ar[r]^-{\iota^G_N} \ar[d]^{\emb} & ([G,G], \dG) \ar[d]^{\emb} \\
0 \ar[r] & \Ker (\bT_{N,G}) \ar[r] & (\bCh(G,N), \bscl_{G,N}) \ar[r]^-{\bT_{N,G}} & (\bCh(G,G), \bscl_G).
}}
\end{align}
 It follows from Proposition~\ref{proposition group lifting} that diagram \eqref{explicit} commutes up to closeness. 

The goal is to show that $\Phi$ is the coarse kernel of $\iota^G_N$. By Corollary~\ref{corollary finite-dimensional}, we have $\dim_ \RR \Ker(\bT_{N,G}) = n$. Thus, if we show that $\sspan_\RR(z'_1, \dots, z'_n) = \Ker(\bT_{N,G})$, then $\Phi'$ is a coarse equivalence, and it follows from Proposition~\ref{proposition 1} that $\Phi$ is the coarse kernel.
Suppose that $\sspan_\RR(z'_1, \dots, z'_n) \ne \Ker(\bT_{N,G})$. By Theorem~\ref{theorem W duality}, there exists $\mu \in \QQQ(N)^G$ whose class in $\WWW(G,N)$ is nonzero and such that $\mu(z'_1) = \dots = \mu(z'_n) = 0$. 
 By the argument in the proof of Proposition~\ref{proposition limit formula} and the quasi-additivity of $\mu$, we have $\mu(z'_i)=h_\mu(\tz_i)$. 
Hence condition $(*)$ implies $\tsigma(\mu) = 0$. However, since $\mu$ does not vanish in $\WWW(G,N)$, Theorem~\ref{theorem core kernel} (or Theorem~\ref{theorem B}) implies $\tsigma(\mu) \ne 0$. This is a contradiction.
\end{proof}

In the rest of this section, we provide several examples of explicit coarse kernels:

\begin{example}
Set
\[ G = \langle g, h \mid [g,h]^2 \rangle \]
and $N = [G,G]$. Then $\dim_\RR \WWW(G,N) = 1$ (see \cite[Remark~4.19]{space}). Thus Theorem~\ref{theorem explicit} and Remark~\ref{remark explicit} imply that
\[ \{ [g,h^m]^2 \mid m \in \ZZ \}\]
is the coarse kernel of $\iota^G_N \colon ([G,N], \dGN) \to ([G,G], \dG)$.  We discussed this example in our Part I paper \cite[Example~5.3]{fa}, without providing a proof.
\end{example}

\begin{example}
Let $l$ be an integer at least $2$. Set
\[ G = \langle a_1, \cdots, a_l, b_1, \cdots, b_l \mid [a_1, b_1] \cdots [a_l, b_l] \rangle\]
and $N = [G,G]$. Then $\dim_\RR \WWW(G,N) = 1$ (see \cite[Theorem~1.1]{space}). Thus Theorem~\ref{theorem explicit} and Remark~\ref{remark explicit} imply that
\[ \{ [a_1, b_1^m] \cdots [a_l, b_l^m] \mid m \in \ZZ \}\]
is the coarse kernel of $\iota^G_N \colon ([G,N], \dGN) \to ([G,G], \dG)$.
\end{example}

\begin{example} \label{example free-by-cyclic}
Let $n$ be an integer at least $2$ and $F_n$ a free group of rank $n$. Here we recall the Andreadakis--Johnson filtration \cite{Andreadakis} of the automorphism group $\Aut(F_n)$ of $F_n$. For every positive integer $q$, the action of $\Aut(F_n)$ on $F_n$ induces the action on $F_n / \gamma_{q+1}(F_n)$. Its kernel is denoted by $\AAA_n(q)$. Then we have the filtration
\[ \Aut(F_n) \geqslant \AAA_n(1) \geqslant \AAA_n(2) \geqslant \cdots \geqslant \AAA_n(q) \geqslant \AAA_n(q+1) \geqslant \cdots.\]
The group $\AAA_n(1)$ is also called the group of IA-automorphisms of $F_n$, and is denoted by $\IA_n$.

Let $\chi \in \Aut(F_n)$, and consider the semi-direct product $G = F_n \rtimes_{\chi}\ZZ$. Let $a_1, \cdots, a_n$ be generators of $F_n$ and $c$ the element of $G$ corresponding to $1 \in \ZZ$. Then, $G$ has the following group presentation
\[
G = F_n\rtimes_{\chi}\ZZ \cong F / R = F(\tilde{a}_1,\ldots,\tilde{a}_{n},\tilde{c}) / \llangle r_1,\ldots,r_{n}\rrangle_F,
\]
where $F = F(\ta_1, \ldots, \ta_n, \tc)$ is the free group generated by $\ta_1, \cdots, \ta_n, \tc$ and
\[ r_i = \tilde{c}\tilde{a}_i\tilde{c}^{-1}\chi(\tilde{a}_i)^{-1}\]
for $i \in \{ 1, \cdots, n\}$. The natural quotient map $\pi \colon F \to G$ sends $\ta_i$ to $a_i$ and $\tc$ to $c$. Here, we regard $\chi \in \Aut(F_n)$ as the automorphism of $F(\ta_1, \cdots, \ta_n) \leqslant F$ as well as the automorphism of $F_n = F(a_1, \cdots, a_n) \leqslant G$.

First we consider the case $\chi \in \AAA_n(1)$ and $\chi$ is atoroidal. Then we have $\ta_i = \chi(\ta_i)$ in $F(\ta_1, \cdots, \ta_n) / \gamma_2(F(\ta_1, \cdots, \ta_n))$, and hence we have $\ta_i = \chi(\ta_i)$ in $F / \gamma_2(F)$. This means that
\[ r_i = \tc \ta_i \tc^{-1} \chi(\ta_i)^{-1} = 1\]
in $F / \gamma_2(F)$, i.e., $r_i \in \gamma_2(F)$ for every $i \in \{ 1, \cdots, n\}$. By \cite[Theorem~4.11]{space}, we have $\dim_\RR \WWW(G, \gamma_2(G)) = n$ in this case. Hence, if we can present $r_i$ as an explicit product of commutators of $F$, then we have an explicit realization of the coarse kernel of $\iota^{G}_{\gamma_2(G)}$.

Now we assume $\chi \in \AAA_n(2)$ in addition. Then we have
\[ \ta_i = \chi(\ta_i) \; \textrm{in $F(\ta_1, \cdots, \ta_n) / \gamma_3(F(\ta_1, \cdots, \ta_n))$,}\] and hence $\ta_i = \chi(\ta_i)$ in $F / \gamma_3(F)$. Hence we have in $F/ \gamma_3(F)$
\[ [\tc, \ta_i] = \tc \ta_i \tc^{-1} \ta_i^{-1} = \tc \ta_i \tc^{-1} \chi(\ta_i)^{-1} = r_i.\]
Therefore, by Theorem~\ref{theorem explicit} and Remark~\ref{remark explicit}, we conclude that
\begin{align*}
&\left\{[c, a_1^{m_1}]\cdots [c,a_n^{m_{n}}]\,\middle|\, (m_1,\ldots,m_{n})\in \ZZ^{n}\right\}\\
={} &\left\{\chi(a_1)^{m_1}a_1^{-m_1}\cdots \chi(a_n)^{m_{n}}a_n^{-m_n}\,\middle|\, (m_1,\ldots,m_{n})\in \ZZ^{n}\right\}
\end{align*}
is the coarse kernel of $\iota^{\Gg}_{\gamma_2(\Gg)}$ for $\Gg=F_n\rtimes_{\chi}\ZZ$.
\end{example}

\begin{example}
Let $\Sigma_g$ be a closed connected oriented surface of genus $g\geq 2$,
and set $\Gamma=\pi_1(\Sigma_g)$.
Write $\Aut_+(\Gamma)$ for the subgroup of $\Aut(\Gamma)$ consisting of
automorphisms whose induced action on $H_1(\Sigma_g;\ZZ)$ preserves
the algebraic intersection form, and set
$\Out_+(\Gamma)=\Aut_+(\Gamma)/\operatorname{Inn}(\Gamma)$,
where $\operatorname{Inn}(\Gamma)$ is the group of inner automorphisms.
By the Dehn--Nielsen--Baer theorem, the action on $\Gamma$ identifies
the orientation-preserving mapping class group $\Mod(\Sigma_g)$
with $\Out_+(\Gamma)$.
The Andreadakis--Johnson filtration of $\Aut(F_n)$ has the following
analogue for $\Mod(\Sigma_g)$.
For $q\in\NN$, let $\Mod_{(q)}(\Sigma_g)$ be the kernel of the natural homomorphism
\[
\Mod(\Sigma_g)\longrightarrow
\Out\bigl(\pi_1(\Sigma_g)/\gamma_{q+1}(\pi_1(\Sigma_g))\bigr).
\] 
This yields the filtration
\begin{equation}\label{eq=Johnson}
\Mod(\Sigma_{g}) \geqslant \Mod_{(1)}(\Sigma_{g})\geqslant \Mod_{(2)}(\Sigma_{g})\geqslant \cdots\geqslant \Mod_{(q)}(\Sigma_{g})\geqslant  \cdots.
\end{equation}
We note that $\Mod_{(1)}(\Sigma_{g})$ equals the Torelli group $\mathcal{I}(\Sigma_{g})$. Also, it follows from the work of Johnson \cite{Johnson} that $\Mod_{(2)}(\Sigma_{g})$ equals the Johnson kernel $\mathcal{K}(\Sigma_{g})$, provided that $g\geq 3$.

Now, we assume that the element $[\chi]$ in $\Mod(\Sigma_{g})$ represented by $\chi\in \Aut_+(\pi_1(\Sigma_{g}))$ is a pseudo-Anosov element. Moreover, we assume that this element $[\chi]$ belongs to the group $\Mod_{(2)}(\Sigma_{g})$ appearing in the filtration \eqref{eq=Johnson}. In particular, the latter assumption implies that $[\chi]\in \mathcal{I}(\Sigma_{g})$. Consider the semi-direct product $\pi_1(\Sigma_{g})\rtimes_{\chi}\ZZ$. 

Consider the following commutative diagram:
\[
\xymatrix{
\Aut(\pi_1(\Sigma_g)) \ar[r]^-{\rho} \ar[d] & \Aut(\pi_1(\Sigma_g) / \gamma_3(\pi_1(\Sigma_g))) \ar[d] \\
\Out(\pi_1(\Sigma_g)) \ar[r] & \Out(\pi_1(\Sigma_g) / \gamma_3(\pi_1(\Sigma_g))).
}
\]
Since $[\chi]$ lies in $\Mod_{(2)}(\Sigma_{g})$, we have $\rho(\chi)\in \Inn (\pi_1(\Sigma_{g})/\gamma_{3}(\pi_1(\Sigma_{g})))$. Take $\sigma\in \pi_1(\Sigma_{g})/\gamma_{3}(\pi_1(\Sigma_{g}))$ such that $\rho(\chi)$ equals the inner conjugation by $\sigma$. Take an arbitrary element $s\in \pi_1(\Sigma_{g})$ that is sent to $\sigma$ by the group quotient map $\pi_1(\Sigma_{g})\twoheadrightarrow\pi_1(\Sigma_{g})/\gamma_{3}(\pi_1(\Sigma_{g}))$. Then, for every $l\in \pi_1(\Sigma_{g})$, we have
\[
\chi(l)\gamma_3(\pi_1(\Sigma_{g}))=sls^{-1}\gamma_3(\pi_1(\Sigma_{g})).
\]
Define $\chi'\in \Aut_+(\pi_{1}(\Sigma_g))$ as the composition of the inner conjugation by $s^{-1}$ and the automorphism $\chi$. Then, $[\chi']=[\chi]$ in $\Mod(\Sigma_g)$, and for every $l\in \pi_1(\Sigma_{g})$, we have
\begin{equation}\label{eq=inngamma_3}
\chi'(l)\gamma_3(\pi_1(\Sigma_{g}))=l\gamma_3(\pi_1(\Sigma_{g})).
\end{equation}

Since $[\chi']=[\chi]$ in $\Mod(\Sigma_g)\cong \mathrm{Out}_+(\pi_1(\Sigma_{g}))$, the group $\pi_1(\Sigma_{g})\rtimes_{\chi'}\ZZ$ is isomorphic to $\pi_1(\Sigma_{g})\rtimes_{\chi}\ZZ$. Hence, we set $G=\pi_1(\Sigma_{g})\rtimes_{\chi'}\ZZ$, and hereafter we treat $G$ instead of $\pi_1(\Sigma_{g})\rtimes_{\chi}\ZZ$. Choose standard generators $a_1,\ldots,a_g,b_1,\ldots,b_g$
of the surface subgroup $\pi_1(\Sigma_g)\leqslant G$; that is,
\[
\pi_1(\Sigma_g)=\langle a_1,\ldots,a_g,b_1,\ldots,b_g \mid [a_1,b_1]\cdots[a_g,b_g]\rangle.
\]
Let $c$ be the element in $G$ corresponding to $1\in \ZZ$. Consider the free group $F=F(\tilde{a}_1,\ldots,\tilde{a}_{g},\tilde{b}_1,\ldots,\tilde{b}_{g},\tilde{c})$. Let $\pi\colon F\twoheadrightarrow G$ be the group quotient map sending $\tilde{a}_1,\ldots,\tilde{b}_{g},\tilde{c}$ to $a_1,\ldots,b_{g},c$, respectively. Then, by \eqref{eq=inngamma_3}, we can take 
\[
\widetilde{\chi}'_{\tilde{a}_1},\ldots,\widetilde{\chi}'_{\tilde{b}_{g}}\in F(\tilde{a}_1,\ldots,\tilde{a}_{g},\tilde{b}_1,\ldots,\tilde{b}_{g})
\]
in such a way that for each $i\in \{1,\ldots,g\}$,
\begin{equation}\label{eq=surfacegamma_3}
\widetilde{\chi}'_{\tilde{a}_i} \gamma_3(F)=\tilde{a}_i\gamma_3(F)\quad \textrm{and}\quad \widetilde{\chi}'_{\tilde{b}_i} \gamma_3(F)=\tilde{b}_i\gamma_3(F)
\end{equation}
as well as
\[
\pi (\widetilde{\chi}'_{\tilde{a}_i})=\chi'(a_i)\quad \textrm{and}\quad \pi (\widetilde{\chi}'_{\tilde{b}_i})=\chi'(b_i).
\]
Then, this group $G=\pi_1(\Sigma_{g})\rtimes_{\chi'}\ZZ$ admits a group presentation
\[
G= F /R =F(\tilde{a}_1,\ldots,\tilde{a}_{g},\tilde{b}_1,\ldots,\tilde{b}_{g},\tilde{c}) / \llangle r_1,r_2,\ldots,r_{2g},r_{2g+1}\rrangle_F,
\]
where $r_i=\tilde{c}\tilde{a}_i\tilde{c}^{-1}(\widetilde{\chi}'_{\tilde{a}_i})^{-1}$ and $r_{g+i}=\tilde{c}\tilde{b}_i\tilde{c}^{-1} (\widetilde{\chi}'_{\tilde{b}_i})^{-1}$ for $i\in \{1,\ldots,g\}$, and $r_{2g +1}=[\tilde{a}_1,\tilde{b}_1]\cdots [\tilde{a}_{g},\tilde{b}_{g}]$. 

By \cite[Theorem~1.2]{space}, we have $\dim_\RR \WWW(G, \gamma_2(G)) = 2g + 1$, and hence we can apply Theorem~\ref{theorem explicit} to this setting. More precisely, by \eqref{eq=surfacegamma_3}, we may choose 
$r_1^{\flat},\ldots,r_{2g+1}^{\flat}$ by setting
\[
r_i^{\flat}=[\tilde{c},\tilde{a}_i]\quad \textrm{and}\quad r_{g+i}^{\flat}=[\tilde{c},\tilde{b}_i]
\]
for every $i\in \{1,\ldots,g\}$, and set $r_{2g +1}^{\flat}=r_{2g +1}$. Therefore, by Theorem~\ref{theorem explicit}  and Remark~\ref{remark explicit},  we conclude that the coarse subspace represented by the set
\[
\left\{\chi'(a_1)^{m_1}a_1^{-m_1}\cdots \chi'(b_{g})^{m_{2g}}b_{g}^{-m_{2g}}[a_1,b_1^{m_{2g+1}}]\cdots [a_{g},b_{g}^{m_{2g+1}}]\,\middle|\, (m_1,\ldots,m_{2g+1})\in\ZZ^{2g+1}\right\}
\]
is the coarse kernel of $\iota^{\Gg}_{\gamma_2(\Gg)}$ for $\Gg=\pi_1(\Sigma_{g})\rtimes_{\chi'}\ZZ$.
\end{example}

\section*{Acknowledgements}

The first-named author, the second-named author, the fourth-named author and the fifth-named author are partially supported by JSPS KAKENHI Grant Numbers JP25K06994, JP24K16921, JP23K12975, and JP26K06799, respectively.
The third-named author is partially supported by JSPS KAKENHI Grant Number JP23K12971 and Research Fellowship Promoting International Collaboration, The Mathematical Society of Japan.

\bibliographystyle{plain} 
\bibliography{reference_scl2}

\end{document}